\documentclass[12pt,reqno]{amsart}
\usepackage{amssymb,delarray}
\usepackage{amsfonts}
\usepackage{epsfig}
\usepackage[all]{xy}
\usepackage{amscd}

\makeindex{}

\newtheorem{thm}{Theorem}
\newtheorem{lem}{Lemma}
\newtheorem{prop}{Proposition}
\newtheorem{claim}{Claim}

\newtheorem{problem}{Problem}

{\catcode`\@=11
\gdef\n@te#1#2{\leavevmode\vadjust{%
 {\setbox\z@\hbox to\z@{\strut#1}%
  \setbox\z@\hbox{\raise\dp\strutbox\box\z@}\ht\z@=\z@\dp\z@=\z@%
  #2\box\z@}}}
\gdef\leftnote#1{\n@te{\hss#1\quad}{}}
\gdef\rightnote#1{\n@te{\quad\kern-\leftskip#1\hss}{\moveright\hsize}}
\gdef\?{\FN@\qumark}
\gdef\qumark{\ifx\next"\DN@"##1"{\leftnote{\rm##1}}\else
 \DN@{\leftnote{\rm??}}\fi{\rm??}\next@}}

\begin{document}

\baselineskip=14.pt plus 2pt 

\title[On complete degenerations of surfaces]{On complete degenerations of
surfaces with ordinary singularities in $\mathbb P^3$}
\author[V.S. Kulikov and  Vik.S.~Kulikov]{V.S. Kulikov and  Vik.S. Kulikov}
\address{Moscow State University of Printing}
\email{vskulikov@mail.ru}

\address{Steklov Mathematical Institute}
\email{kulikov@mi.ras.ru}

\thanks{This
research was partially supported by grants of NSh-9969.2006.1, RFBR
08-01-00095 and RFBR 06-01-72017-$\text{MNTI}_{-}a$. \newline The
research was started during the stay of the second author in Centro
di Ricerca Matematica Ennio de Giorgi (Program: Groups in Algebraic
Geometry).}

\keywords{}

\begin{abstract}
We investigate the problem of existence of degenerations of surfaces
in $\mathbb P^3$ with ordinary singularities into plane arrangements
in general position.
\end{abstract}

\maketitle

\setcounter{tocdepth}{2}


\def\st{{\sf st}}

\section*{Introduction.}
In the article we investigate degenerations of surfaces in
$\mathbb P^3$ with ordinary singularities. To begin, consider the
classical prototype of this situation, namely, degenerations of
plane algebraic curves. As is known, any smooth projective curve
can be projected to $\mathbb P^2$ onto a nodal curve $C$ --- a
curve with ordinary double points
--- nodes (singularities of type $A_{1}$). According to Severi
Theorem (\cite{Sev},\cite{H}) any nodal plane curve of degree $m$
can be degenerated into an arrangement of $m$ lines in general
position  $\mathcal L=\cup_{i=1}^m L_i\subset \mathbb P^2$. Such
degeneration defines a subset $S_{0}$ of the set  $S$ of double
points of the curve $\mathcal L$, which consists of the limit
double points of the curve $C$. Conversely, for any subset $S_{0}
\subset S$ there exists a degeneration of nodal curves (may be,
reducible) for which  $S_{0}$ is the set of limit double points
(\cite{GLS}), in other words, there exists a smoothing of double
points  $S \setminus S_{0}$ of the curve $\mathcal L$ (preserving
double points lying in $S_{0}$).

In dimension 2, the analog of nodal curves are surfaces in
${\mathbb{P}}^{3}$ with ordinary singularities. Only such
singularities appear under generic projections of a smooth surface
$X \subset \mathbb P^r$ to ${\mathbb{P}}^{3}$ (see, for example,
\cite{Mo} and \cite{G-H}). Let $Y \subset{\mathbb{P}}^{3}$ be a
surface of degree $m$ with ordinary singularities. This means that
the singular set ${\rm Sing}Y$  is the double curve $D$; the curve
$D$ itself is smooth except the finite set of triple points $T$; the
generic point $y\in D$  is nodal on $Y$ (locally defined by equation
$xy=0$); the points $y\in T$ are triple points of $Y$ (locally,
$xyz=0$); besides, $Y$ has a finite set of pinch points at which $Y$
is locally defined by equation $x^{2}=y^{2}z$.

Maximally degenerate (reducible) surface of degree $m$ with
ordinary singularities is an arrangement  $\mathcal P=P_1 \cup
\dots \cup P_m \subset{\mathbb{P}}^{3}$ of $m$ planes in general
position. In this case the double curve ${\rm Sing
}{\,\mathcal{P}}={\mathcal{L}}$ consists of ${m\choose 2}$ lines
$L_{i,j}= P_i\cap P_j$,
$$
{\mathcal{L}} =\bigcup_{1\leq i<j\leq m}L_{i,j} ,
$$
on which there lie  ${m\choose 3}$  triple points,
$$
 {\mathcal{T}} =\cup T_{i,j,k}, \quad T_{i,j,k}=P_i\cap P_j\cap P_k,\quad 1\leq i<j<k\leq m .
$$

A degeneration of surfaces (or, in general, of varieties) of given
type  is a one-dimensional family , or for  brevity, its zero
fibre, generic fibre of which is a surface of given type (usually,
smooth). We consider families, all fibres of which are surface in
$\mathbb P^3$ with ordinary singularities. A degeneration is a
surface, which has "more"\, singularities than generic fibre. More
precisely, a {\it degeneration } of a surface $Y \subset
{\mathbb{P}}^{3}$ with ordinary singularities is a flat family of
embedded surfaces $Y_{u} \subset {\mathbb{P}}^{3}$ with ordinary
singularities, parametrized by points  $u\in U$ of a smooth curve
$U$ (or even, of a disk $U\subset \mathbb C$), such that

\begin{itemize} \item[$(i)$] for the generic point  $u\in U$ the
fibres $Y_{u}$ have singularities of the same type as $Y=Y_{u_1}$,
and the fibre  $Y_{u_{0}}$  over the point $u_{0}=0$ is called {\it
degenerate};
\item[$(ii)$] there is a flat family $D_{u}\subset Y_{u}$, where
$D_{u}$ is the double curve of the surface $Y_{u}$ for $u\neq
u_{0}$,  and the curve $\mathcal D = D_{u_{0}} \subset Y_{u_0}$ is
called the {\it limit  double curve};
\item[$(iii)$]  there is a flat family  $T_{u}$, $u\in U$, where $T_{u}\subset
D_{u}$ is the set of triple points  of the surface  $Y_{u}$ for
$u\neq u_{0}$, and $T_{u_{0}}={\mathcal T}_{3}$ is the set of
triple points  of the curve $\mathcal D $.
\end{itemize}
A degeneration is called  {\it complete} if the degenerate fibre
$Y_{u_{0}} = \mathcal P$ is a plane arrangement in general position.

A complete degeneration defines a limit double curve  $\mathcal D
\subset\mathcal L$.  Conversely, if a line arrangement $\mathcal D
\subset\mathcal L$ is selected and there is a complete degeneration
of surfaces $Y_{u} \subset {\mathbb{P}}^{3}$ (not necessary
irreducible) with ordinary singularities such that
$D_{u_{0}}=\mathcal D $, then we say that the plane arrangement
$\mathcal P$ {\it is smoothed outside of} $\mathcal D $.

As in the case of curves, there are two questions. {\it Whether
every surface $Y \subset{\mathbb{P}}^{3}$ with ordinary
singularities can be completely degenerated?}  {\it Whether every
line arrangement $\mathcal D \subset\mathcal L$ of  a plane
arrangement $\mathcal P$ can be smoothed outside of $\mathcal D
$}?

The problem of degeneracy of smooth surfaces $X\subset \mathbb
P^N$ into a plane arrangements was investigated previously, but in
another setting, which is analogous to the Zeuthen's problem for
curves (the degenerate surfaces are surfaces of Zappa, or so
called Zappatics; a survey of obtained results about Zappatic
surfaces can be found in \cite{Cil}, see also \cite{Cil-M}).

In spite of the fact that in many cases surfaces with ordinary
singularities can be completely degenerated (see $\S$ \ref{deg4}
and $\S$ \ref{cmlt}), we shall show in the article that in general
the answers to the questions formulated above are in negative. The
reason for impossibility to completely degenerate a surface with
ordinary singularities can lie  as in the fact that it is
impossible to degenerate the double curve of the surface to  a
line arrangement $\mathcal D$ (see $\S$ \ref{examp}), just as in
the fact that the corresponding limit pair $(\mathcal P,\mathcal
D)$ does not exist, although the double curve of the surface can
be degenerated to a line arrangement (see $\S$ \ref{example2}). In
addition we show that there exist pairs $(\mathcal P_m,\mathcal
D)$, which can not be smoothed outside of a fixed line arrangement
$\mathcal D $. In some cases such pairs can not be smoothed for
all $m=\deg \mathcal P_m$ (see $\S$ \ref{pot2}), and in the other
cases such pairs can  be smoothed only for sufficiently large $m$
(see $\S$ \ref{pot1}).

\section{Expression of numerical characteristics of a surface in terms of degeneration}
There are two classical ways to study smooth surfaces $X \subset
\mathbb P^r$: to consider its generic projection to $\mathbb P^3$
or onto $\mathbb P^2$ respectively. We consider a few more general
situation. We begin with a surface $Y \subset{\mathbb{P}}^{3}$
with ordinary singularities, and the normalization $\mathfrak{n}:
X\rightarrow Y$ of $Y$ and its composition with projection of $Y$
onto $\mathbb P^2$ define a generic covering $X \rightarrow
\mathbb P^2$. The definition of a complete degeneration makes it
possible to express numerical characteristics of the surface $Y$
in terms of its degeneration. Invariants of a surface $X$ can be
expressed as by means of  numerical characteristics of the surface
$Y$, just as  by means of  numerical characteristics of a generic
covering $X \rightarrow \mathbb P^2$. This gives an expression of
numerical characteristics of the covering in terms of
degeneration.

\subsection{Numerical characteristics of a surface with ordinary
singularities.} \label{type1} Let $Y \subset{\mathbb{P}}^{3}$ be
an irreducible surface with ordinary singularities,  $D=D_{1}\cup
\dots \cup D_{k}$ the double curve $g_{i}=g(D_{i})$ its geometric
genus, $d_{i}= {\rm deg}D_{i}$ the degree of an irreducible
component  $D_{i}$, $T$ the set of triple points,  $\Omega$ the
set of pinch points. The main numerical characteristics of the
imbedding $Y \subset{\mathbb{P}}^{3}$ are:

$m={\rm deg\,}Y$ -- the degree of the surface,

$k$ -- the number of irreducible components of the double curve,

$\bar g =\sum_{i=1}^{k}g_{i}$ -- the geometric genus of the double
curve,

$\bar d =\sum_{i=1}^{k}d_{i}$ -- the degree of the double curve,

$t=\sharp(T)$ -- the number of triple points,

$\omega =\sharp (\Omega)$ -- the number of pinch points.

The collection of numerical data  ${\rm type\,}(Y) =(m,\bar d, k,
\bar g, t)$ we call {\it the type of the surface} $Y$.

Let $\mathfrak{n} : X\rightarrow Y$ be a normalization of the
surface  $Y$. The surface $X$ is smooth. Its invariants  are
expressed in terms of numerical characteristics of the surface $Y$
by the following formulae (see {\cite{G-H}}):

the intersection number of the canonical class is equal to
\begin{equation} \label{K21}
K^2_X=m(m-4)^2-(5m-24)\bar d +4(\bar g -k)+9t ,
\end{equation}

the topological Euler number  $e(X)=c_{2}(X)$ is equal to
\begin{equation} \label{E1}
e(X)=m^2(m-4)+6m-(7m-24)\bar d +8(\bar g -k)+15t ,
\end{equation}
From these formulae and the Noether's formula we obtain the Euler
characteristic  $\chi ({\mathcal{O}}_{X})$:
\begin{equation} \label{eul} \chi ({\mathcal{O}}_{X})=\frac{1}{6}m(m^2-6m+11)-(m-4)\bar d+(\bar g
-k)+2t.\end{equation}

The number of pinches  is equal to  (see {\cite{G-H}}):
\begin{equation} \label{w}
\omega =2\bar d (m-4) -6t -4(\bar g -k) .
\end{equation}

As is known (\cite{B-G}), the arithmetic genus of a reduced curve
$C$ is equal to
\begin{equation} \label{pa}
p_{a}(C) =g+\delta -r+1 ,
\end{equation}
where $g$ is its geometric genus,  $\delta$ is the sum of
$\delta$-invariants of the singularities,  $r$ is the number of
irreducible components. For the double curve $D$, which has only
triple points  ($\delta =3$),  formula (\ref{pa}) gives
\begin{equation} \label{paD}
p_{a}(D)=\bar g +3t- k +1.
\end{equation}

\subsection{Numerical data for description of a curve  $\mathcal D \subset\mathcal
L$} \label{type2} Let $\mathcal P \subset{\mathbb{P}}^{3}$ be an
arrangement of $m$ planes in general position, ${\mathcal{L}}
=\cup_{1\leq i<j\leq m}L_{i,j}$ its double curve. Let a curve
$\mathcal D\subset \mathcal L$ be the union of $\bar d$ lines
$L_{i,j}$, and $\mathcal R$ the union of $d$ remainder lines (the
curves $\mathcal D$ and  $\mathcal R$ also are called line
arrangements),
\begin{equation} \label{L}
\mathcal L =\mathcal D \cup \mathcal R, \quad \deg \mathcal D
=\bar d,\quad \deg \mathcal R = d, \quad d+\bar d ={m\choose 2}.
\end{equation}
With respect to the partition $\mathcal L=\mathcal D \cup \mathcal
R$ the triple points of the curve  $\mathcal L$ are decomposed
into 4 types:
$$\mathcal T=\mathcal T_3\sqcup\mathcal T_2\sqcup\mathcal
T_1\sqcup\mathcal T_0,$$ where $\mathcal T_3$ consists of those
points, which are triple on the curve  $\mathcal D$; $\mathcal T_2$
consists of those points, which are double on the curve $\mathcal
D$; $\mathcal T_1$  consists of those points, which are non-singular
on the curve $\mathcal D$ (and are double on the curve $\mathcal
R$); $\mathcal T_0$ consists of those points, which do not lie on
$\mathcal D$ (and are  triple on the curve $\mathcal R$). Denote by
$\tau_{3}$, $\tau_{2}$, $\tau_{1}$, and $\tau_{0}$ the number of
points in corresponding sets. We have
\begin{equation} \label{tau}
\tau =\tau_{3}+\tau_{2}+\tau_{1}+\tau_{0}= {m\choose 3} .
\end{equation}
By the formula for arithmetic genus (\ref{pa}), we have
\begin{equation} \label{paD0}
p_{a}(\mathcal D) =\tau_{2}+3\tau_{3}-\bar d +1 .
\end{equation}
(the invariant  $\delta =1$ for double points, and $\delta =3$ for
triple points).

The numbers $\tau_{i}$ are related as follows. On each line
$L_{i,j}= P_i\cap P_j$ there are $m-2$ triple points
--- the intersection points  with planes  $P_{k}$, $k\neq i,j$.
Let us sum the numbers of triple points of $\mathcal L$ lying on
$\bar d$ lines $L_{i,j}\subset \mathcal D$. On one hand, we get
$(m-2)\bar d$. On the other hand, we get
$3\tau_{3}+2\tau_{2}+\tau_{1}$, because the triple points of
$\mathcal D$ are counted 3 times, the double points
--- 2, and non-singular --- 1. We get
\begin{equation} \label{tau123}
\tau_{1} +2\tau_{2}+3\tau_{3} = (m-2)\bar d .
\end{equation}
the analogous calculation for the curve $\mathcal R$ gives
\begin{equation} \label{tau210}
\tau_{2} +2\tau_{1}+3\tau_{0} = (m-2)d .
\end{equation}

The data collection ${\rm type\,}(\mathcal P ,\mathcal D) =(m,\bar
d, k,\tau_{2},\tau_{3})$ is called  {\it the type of pair }
$(\mathcal P ,\mathcal D)$. Here $k$ is the number of connected
components of the curve $\mathcal D \setminus {\mathcal{T}}_{3}$,
where ${\mathcal{T}}_{3}$ is the set of triple points of $\mathcal
D$.

\subsection{The graph of a line arrangement}
As is known, one can associate to an algebraic variety, and in
particular, to a divisor with normal crossings, a polyhedron by the
rule: the vertices correspond to the irreducible components; two
vertices are connected by an edge if the components have a non empty
intersection; three vertices span a triangle if the corresponding
components have a non empty intersection and etc. In particular, to
an arrangement  $\mathcal P$ of $m$  planes in general position in
$\mathbb{P}^{3}$, we associate a polyhedron, which is the
two-dimensional skeleton of a $(m-1)$-dimensional simplex in
$\mathbb{R}^{m-1}$ (or a standard $(m-1)$-simplex in
$\mathbb{R}^{m}$). Denote by $\Gamma (\mathcal L)$ a graph, which is
the one-dimensional skeleton of this polyhedron.

To any line arrangement $\mathcal D \subset \mathcal L$ we
associate a graph $\Gamma (\mathcal D)$ --- {\it the graph of the
curve} $\mathcal D$, which is a  subgraph of $\Gamma (\mathcal
L)$, consisting of the union of edges corresponding to the lines
$L_{i,j} \subset \mathcal D$. If $\mathcal D \subset \mathcal L$
is a line arrangement, then we denote by $\mathcal R$ the
complementary line arrangement in $ \mathcal L$.

We say that a graph  $\Gamma$ is {\it realizable}, if it satisfies
the following conditions: it doesn't contain isolated vertices, it
doesn't  contain  simple loops (that is, edges with the same
source and endpoint), and any two vertices are connected by at
most one edge. It is obvious,  a graph $\Gamma$ is the graph  of a
set of double lines $\mathcal D$ of a plane arrangement $\mathcal
P$ in general position  if and only if $\Gamma$ is realizable and
the number of its vertices is not more than $\deg \mathcal P$.

The graph  $\Gamma (\mathcal D)$ codes  numerical characteristics
of $\mathcal D$. The number of vertices of $\Gamma (\mathcal D)$
is equal to the number $\bar d$ of lines composing the curve
$\mathcal D$. The number  $\tau_{3}$ of triple points of $\mathcal
D$, obviously,   is equal to the number of triangles of $\Gamma
(\mathcal D)$ (by definition, a triangle in a graph is three
vertices and three edges connecting these vertices).

Denote by $v(P)$ the valence of a vertex $P$ of $\Gamma (\mathcal
D)$, that is,  the number of  edges outgoing of $P$.
\begin{lem}\label{l-v}
The number of the double points of the curve $\mathcal D$ is equal
to
\begin{equation} \label{v}
\tau_{2} =\sum_{P\in \Gamma (\mathcal D)}\frac{1}{2}v(P)(v(P)-1)
-3\tau_{3}
\end{equation}
\end{lem}
\proof A vertex  $P$ of valence $v(P)$ corresponds to a plane, on
which $v(P)$ lines lie. These lines intersect at
$\frac{1}{2}v(P)(v(P)-1)$ double points. The sum of numbers of
these double points is  the number of double points of $\mathcal
D$, if $\mathcal D$ has no triple points. If $\mathcal D$ has
triple points, then each triple point $T_{i,j,k}= P_{i}\cap
P_{j}\cap P_{k}$, being a  double point on each of three planes
$P_{i}$, $P_{j}$, $P_{k}$, contribute 3 to the first summand of
the formula (\ref{v}). To calculate $\tau_{2}$, we have to remove
this contribution, that is, to subtract $3\tau_{3}$. \qed

A graph $\bar{\Gamma}(\mathcal R)$, consisting of all  vertices of
the graph $\Gamma(\mathcal L)$ and all edges of the graph
$\Gamma(\mathcal R)$ is called {\it an augmentation of}
$\Gamma(\mathcal R)$ ({\it with respect to} $\mathcal L$). A graph
$\Gamma(\mathcal R)$ is called {\it augmented}, if
$\Gamma(\mathcal R) = \bar\Gamma(\mathcal R)$.

Let $\mathcal D \subset \mathcal L$ be the limit double curve of a
complete degeneration of a surface $Y$ with ordinary singularities
in $\mathbb{P}^{3}$. It is not difficult to prove the following
lemma.

\begin{lem}\label{R}
The number of irreducible components of a surface $Y$  with
ordinary singularities in $\mathbb{P}^{3}$ is equal to the number
of connected components of the graph $\bar\Gamma(\mathcal R)$. In
particular, a surface $Y$ is irreducible, if the graph
$\bar\Gamma(\mathcal R)$ is connected.
\end{lem}

In connection with this lemma, a line arrangement $\mathcal D$,
more precisely, a pair  $(\mathcal L , \mathcal D)$ is called {\it
irreducible}, if the graph $\bar\Gamma(\mathcal R)$ is connected.

In the case of a connected graph we use a symbol
$$
\Gamma_{v_{1},v_{2},\dots}^{\bar d,\tau_{3}}
$$
for the type of the graph $\Gamma(\mathcal D)$ (or for the line
arrangement $\mathcal D$ itself). Here  $v_{1},v_{2},\dots$ are the
numbers of vertices of  valence 1,2 and etc., $\bar d$ is the degree
of $\mathcal D$, that is, the number of edges of the graph,
$\tau_{3}$ is the numbers of triple points of $\mathcal D$, that is,
the number of triangles in $\Gamma(\mathcal D)$. If the graph
$\Gamma(\mathcal D)$ is not connected, then we use the analogous
notation for connected components and separate data for each
component by parentheses. If several components have the same data,
we use a multiplicative way of writing. For example, if $\mathcal D$
consists of three connected components, two of which are chains of
two lines, and the third consists of four lines, three of which
intersect in a triple point and the fourth line intersect two of the
three lines, then the graph $\Gamma(\mathcal D)$ has type
$\Gamma^{(2,0)^2(4,1)}_{(2,1)^2(1,2,1)}$ (see Fig.1).

\begin{figure}[h]
\begin{picture}(0,100)(160,-90)
\put(100,-12){\circle*{3}} \put(100,-11){\line(1,0){30}}
\put(130,-12){\circle*{3}} \put(130,-11){\line(1,0){30}}
\put(160,-12){\circle*{3}} \put(100,-23){\circle*{3}}
\put(100,-22){\line(1,0){30}} \put(130,-23){\circle*{3}}
\put(130,-22){\line(1,0){30}} \put(160,-23){\circle*{3}}

\put(175,-2){\line(0,-1){30}} \put(175,-2){\circle*{3}}
\put(175,-32){\circle*{3}} \put(190,-17){\circle*{3}}
\put(175,-2){\line(1,-1){15}} \put(175,-32){\line(1,1){15}}
\put(190,-16){\line(1,0){25}} \put(215,-17){\circle*{3}}
\put(140,-60){$\Gamma^{(2,0)^2(4,1)}_{(2,1)^2(1,2,1)}$}

\put(150,-85){$\text{\rm Fig}.\, 1$}
\end{picture}
\end{figure}

\subsection{Numerical characteristics of the limit double curve $\mathcal D \subset\mathcal L$}
The connection between numerical characteristics of double curves
$\mathcal D $ and $D$ of a complete degeneration $Y_u$, $u\in U$,
is given by the following
\begin{prop} \label{prop1}
There are following equalities for the limit double curve
$\mathcal D \subset\mathcal L$ of a complete degeneration $Y_{u}$,
$u\in U$, of surfaces with ordinary singularities:

\begin{equation} \label{degD}
\bar d ={\rm deg}\,D  ;
\end{equation}
\begin{equation} \label{tau3}
\tau_{3} = t ;
\end{equation}
\begin{equation} \label{tau2}
\tau_{2} =  \bar d + \bar g -k ;
\end{equation}
\begin{equation} \label{wtau}
\omega =2\tau_{1} ;
\end{equation}
\begin{equation} \label{kŒ}
k = \sharp\, C(\mathcal D\setminus \mathcal T_3) ,
\end{equation}
where $\sharp\, C(\mathcal D\setminus \mathcal T_3)$ is the number
of connected components of the curve $\mathcal D \setminus
{\mathcal T}_{3}$.

Thus, if a surface $Y=Y_{u_1}$ is irreducible, then the types:
${\rm type\,}(Y)$ and ${\rm type\,}(\mathcal P ,\mathcal D)$,
uniquely restore each other.
\end{prop}
\proof Formula  (\ref{degD}) is valid because the degree $\deg
D_{u}$ is constant for a flat family $D_{u}$. Analogously, formula
(\ref{tau3}) follows from the flatness of the  family $T_{u}$.
Formula  (\ref{tau2}) follows from the constancy of the arithmetic
genus $p_{a}(D_{u})$ in  a flat family $D_{u}$ and from  formulae
(\ref{paD}) and (\ref{paD0}). To prove (\ref{wtau}) substitute the
expression for $( \bar g -k )$ from (\ref{tau2}) to formula
(\ref{w}). We have
$$
\omega = 2 \bar d (m-4) -6\tau_{3} +4( \bar g -\tau_{2} ) =2 \bar
d (m-2) -4\tau_{2} -6\tau_{3} ,
$$
and from (\ref{tau123}) we obtain that $\omega =2\tau_{1}$.

Formula (\ref{kŒ})  follows from flatness of the family of curves
$D_{u} \setminus T_{u}$ and from constancy of the number of
connected components in such a family.  \qed

The number $k$ of irreducible components of the double curve
$D\subset Y$ (formula (\ref{kŒ}))) can be expressed in terms of
the graph  $\Gamma(\mathcal D)$ of the limit curve $\mathcal D$ as
follows. By definition, a path in $\Gamma(\mathcal D)$ is a
sequence of edges, for which the initial vertex of the next edge
is the the endpoint of the previous  edge. A path is called {\it
prohibited}, if two of its successive edges are sides of a
triangle in $\Gamma(\mathcal D)$. Call two edges of the graph
$\Gamma(\mathcal D)$ {\it equivalent}, if they are  edges of some
non prohibited path, and continue this relation to equivalence
relation on the set edges of the graph $\Gamma(\mathcal D)$.
Immediately one obtains the following lemma.
\begin{lem}\label{k}
The number $k$ of connected components of the curve $D$ is equal
to the number of equivalence classes of edges of the graph
$\Gamma(\mathcal D)$.
\end{lem}

The formulae (\ref{tau3}) and (\ref{tau2}) permit to express the
invariants of the surface $X$ in terms of numerical
characteristics of the pair $(\mathcal P,\mathcal D)$:
\begin{equation} \label{K2}
K^2_X=m(m-4)^2-(5m-20)\overline d+4\tau_2+9\tau_{3},
\end{equation}
\begin{equation} \label{E}
e(X)=m^2(m-4)+6m-(7m-16)\overline d+8\tau_2+15\tau_{3},
\end{equation}
or, taking into account (\ref{tau123}),
\begin{equation} \label{K2tau}
K^2_X=m(m-4)^2 +10\bar d -5\tau_{1} -6\tau_2 -6\tau_{3},
\end{equation}
\begin{equation} \label{Etau}
e(X)=m^2(m-4)+6m +2\bar d -7\tau_{1} -6\tau_2 -6\tau_{3}
\end{equation}
It follows from the Noether's formula and formula (\ref{tau}) that
\begin{equation} \label{pa1} \chi ({\mathcal{O}}_{X})=\bar d+\tau_{0} -\frac{1}{2}m(m-1) .\end{equation}

\section{Generic projections onto plane}
\subsection{The discriminant of projection}
Let $Y\subset \mathbb P^3$ be a surface (not necessary irreducible)
of degree $m$ with ordinary singularities  and let
$\text{pr}:\mathbb P^3\to\mathbb P^2$ be a linear projection with
the center at a point $o\not\in Y$. Choose coordinates in
$\mathbb{P}^{3}$ such that $o =(0:0:0:1)$. Then the projection
$\text{pr}$ is given by
$$(x_{1}:x_{2}:x_{3}:x_{4}) \mapsto
(x_{1}:x_{2}:x_{3}) \in \mathbb{P}^{2}$$ and $Y$ has an equation
$$
h(x_4)= x_4^{m}+\sum_{j=0}^{m-1}a_{j}(x_1,x_2,x_3)x_4^{j} =0.
$$
The discriminant $\Delta(x_1,x_2,x_3)$ of the polynomial $h(x_4)$,
as  a polynomial in $x_4$, is a homogeneous polynomial in
variables $x_1,x_2,x_3$ of degree $m(m-1)$, and defines in
$\mathbb P^2$ {\it a discriminant divisor} $\Delta$ of
$p=\text{pr}_{\mid Y}$, given by equation $\Delta(x_1,x_2,x_3)=0$.

\begin{prop} \label{dis} Let a projection ${\rm pr}:\mathbb P^3\to \mathbb P^2$ be such that
\begin{itemize} \item[$(i)$] the composition $f=p\circ {\mathfrak{n}}
: X\rightarrow {\mathbb{P}}^{2}$ is unramified over generic points
of components of the image $\bar D=\text{pr}(D)$ of the double
curve $D\subset Y$, where ${\mathfrak{n}} : X\rightarrow Y$ is a
normalization;
\item[$(ii)$] the ramification index of $f$ at generic points of
its ramification curve $R\subset X$ is equal to two;
\item[$(iii)$] the restriction of $f$ to the curve $R\subset X$ is of degree one.
\end{itemize}
Then the polynomial $\Delta(x_1,x_2,x_3)$ factors as:
$$\Delta(x_1,x_2,x_3)=\beta(x_1,x_2,x_3)\cdot \rho(x_1,x_2,x_3)^2,$$ where polynomials
$\beta(x_1,x_2,x_3)$ and $\rho(x_1,x_2,x_3)$ have no multiple
factors,  the equation of the branch curve $f(R)=B\subset \mathbb
P^2$ of $f$ is $\beta (x_1,x_2,x_3)=0$, and $\rho (x_1,x_2,x_3)=0$
is the equation of the curve $\bar D$.
\end{prop}
\proof The surface $Y$ is covered by three charts $U_{i}=\{ x_{i}
\neq 0\}$, $i=1,2,3$, $U_{i}=\mathbb{C}^{3}$. In each of these
charts the projection $\text{pr}$ is of the form $(x,y,z) \mapsto
(x,y)$. For example, in chart $U_{3}=\{ x_{3} \neq 0\}$ one takes
$\frac{x_{1}}{x_{3}} =x$, $\frac{x_{2}}{x_{3}} = y$,
$\frac{x_{4}}{x_{3}} = z$. In this chart the surface $Y$ has
equation
$$
h(z)= z^{m}+\sum_{j=0}^{m-1}\bar a_{j}(x,y)z^{j} =0,
$$
and the equation of the discriminant divisor is $\Delta(x,y,1)=0$.

It follows from conditions on  the projection $\text{pr}$ that for
each point $q$, lying in one of the three charts isomorphic to
$\mathbb C^2$ (except for a finite set of points --- the singular
points of $B$ and the images of pinches and triple points of $Y$),
there are an analytic neighborhood $V\subset \mathbb C^2$ of this
point and such analytic coordinates $u,v$ that the polynomial
$h(z)$ is a product of $m$ factors of the form $z-\alpha_j(u,v)$
if $q\not\in B$, or of $m-2$ factors of the same form and one of
the form $(z-g(u,v))^2-u$ if $q\in B$, where $\alpha_j(u,v)$,
$g(u,v)$ are some analytic  functions  and $u=0$ is a local
equation of the curve $B$. Furthermore, as is known, the
discriminant of a polynomial $h(z)=\prod_{j=1}^m(z-\alpha_j)$ is
equal to $\prod_{i<j}(\alpha_i-\alpha_j)^2$. Consequently, the
discriminant $\Delta(x,y,1)$ vanishes only at the points of $B$
and $\bar D=\text{pr}(D)$. Moreover, as $Y$ has only ordinary
singularities, it follows from the conditions on the projection
$\text{pr}$ that the equations of the irreducible components of
the curve $\bar D$ appear in $\Delta(x,y,1)$ with multiplicity
two, and in the equations of the irreducible components of $B$
appear with multiplicity one.  \qed \\

In particular, if $Y=\mathcal P$ is a plane arrangement, then the
discriminant divisor of restriction of the projection $\text{pr}$
to $\mathcal P$ equals $\Delta=2\text{pr}(\mathcal
L)=2\bar{\mathcal L}$. Moreover, if the projection $\text{pr}$
satisfies the conditions of Proposition \ref{dis} for each surface
$Y_u$, $u\in U$, of a complete degeneration of surfaces with
ordinary singularities,  then we obtain three flat families of
curves: $\bar D_u=\text{pr}(D_u)$, $B_u$ and $\Delta_u=B_u+2\bar
D_u$, where $\bar D_{u_0}=\text{pr}(\mathcal D)=\bar{\mathcal D}$
and $B_{u_0}=2\text{pr}(\mathcal R)=2\bar{\mathcal R}$.

\subsection{Generic coverings of the plane and generic projections of surfaces with ordinary
singularities} \label{1.4} Let $X$ be a smooth algebraic surface,
not necessary irreducible. Recall  (\cite{Ku}) that a finite
covering $f: X\rightarrow {\mathbb{P}}^{2}$ is called {\it
generic}, if it is like a generic projection of a projective
surface onto the plane, that is, satisfies the following
properties
\begin{itemize} \item[$(i)$] the branch curve $B\subset \mathbb P^2$
is cuspidal, that is, has ordinary cusps and nodes, as the only
singularities;
\item[$(ii)$]  $f^{*}(B)=2R+C$, where  the ramification curve $R$
is smooth, and the curve $C$ is reduced;
\item[$(iii)$]  $f_{\mid R}:R\to B$ is the normalization of $B$.
\end{itemize}

\begin{prop} \label{prop2} Let $X$ be  a  smooth irreducible projective surface.
Then the branch curve $B\subset \mathbb P^2$  of a generic
covering $f:X\to\mathbb P^2$ of degree $m\geq 2$ is irreducible.
\end{prop}
\proof The statement is obvious if $\deg f=2$.

Let $\deg f\geq 3$. A generic covering $f:X\to\mathbb P^2$ branched
along a curve $B$ defines (and is defined by) an epimorphism $\mu
:\pi_1(\mathbb P^2\setminus B) \rightarrow \mathcal S_m$ to the
symmetric group $\mathcal S_m$ such that the image  $\mu (\gamma)$
of each geometric generator $\gamma\in \pi_1(\mathbb P^2\setminus
B)$ (that is, of each simple circuit around the curve $B$)  is a
transposition in $\mathcal S_m$ (see, for example, \cite{Ku}). If
$B$  splits  into the union of two curves $B_1$ and $B_2$, then
$B_1$ and $B_2$ meet each other transversally at non-singular
points, since $B$ has only nodes and ordinary cusps, as the only
singularities.  Therefore (see, for example, \cite{Ku1}) the
elements of the group $\Gamma _1$, generated by the simple circuits
around  $B_1$, commute with the elements of the group $\Gamma _2$,
generated by the simple circuits around $B_2$. On the other hand, it
is easy to see that the elements of two nontrivial subgroups
$\mu(\Gamma_1)\subset \mathcal S_m$ and $\mu(\Gamma_2)\subset
\mathcal S_m$, generated by transpositions, can not generate
$\mathcal S_m$, if the elements of the group $\mu(\Gamma _1)$
commute with the elements of the group $\mu(\Gamma _2)$.
\qed \\

As a corollary of Proposition \ref{prop2} we obtain that the
number of irreducible components of the branch curve $B$ of a
generic covering $f:X\to \mathbb P^2$ is equal to the number of
irreducible components $X_i$ of $X$ such that $\deg f_{\mid
X_i}\geq 2$.

The main numerical characteristics of a generic covering  $f$ of
projective plane by an irreducible surface are:

$m={\rm deg\,}f$ -- the degree of the covering,

$g=g(B)$ -- the geometric genus of the branch curve,

$2d ={\rm deg\,}B$ -- the degree of the branch curve,

$n$ -- the number of nodes of $B$,

$c$ -- the number of cusps.

The invariants of a surface $X$ are expressed in terms of
numerical characteristics of a generic covering  by the following
formulae (see \cite{Ku}):
\begin{equation} \label{k2}
K^2_X=9m-9d+g-1,
\end{equation}
\begin{equation} \label{e}
e(X)=3m+2(g-1)-c,
\end{equation}

Now let $Y \subset{\mathbb{P}}^{3}$ be a surface with ordinary
singularities and ${\mathfrak{n}} : X\rightarrow Y$ be its
normalization. A projection ${\rm pr}:\mathbb P^3\to \mathbb P^2$
is called {\it generic with respect to} $Y$, if it is generic for
the double curve $D\subset Y$ (and, in particular, ${\rm pr}(D)$
has only  $t$ triple points and some nodes as singularities) and
the composition $f=p\circ {\mathfrak{n}} : X\rightarrow
{\mathbb{P}}^{2}$ is a generic covering of the plane, where
$p={\rm pr}_{\mid Y}$.

\begin{prop} \label{prop3}  Let $Y_{u}$, $u\in U$, $\dim U=1$, be a complete degeneration
of surfaces of degree $m$ with ordinary singularities. Then for
for almost all generic projections ${\rm pr}:\mathbb P^3\to
\mathbb P^2$ of the plane arrangement $\mathcal P=Y_{u_0}$ the
projections ${\rm pr}$ are  generic with respect to $Y_u$ for all
$u\in U$, may be, except a finite number of values $u$.
\end{prop}
\proof According to definition of a degeneration the surfaces
$Y_{u}$, $u\in U$, are the fibres of a restriction to $\mathcal
Y\subset \mathbb P^3\times U$ of the projection onto the second
factor.

By definition of a generic projection, applied to a plane
arrangement, it follows that there exist a finite covering $\cup
U_i=\mathbb P^2$  by small open balls $U_i$ with centers at points
$p_i$  and an open neighborhood $V$ of $\mathcal P$  such that the
number of connected components of the intersection ${\rm
pr}^{-1}(U_i)\cap V$  equals respectively to: $m$, if $p_i\not\in
\bar{\mathcal L}={\rm pr}(\mathcal L)$, $m-1$, if $p_i$ is a
non-singular point of $\bar{\mathcal L}$, and $m-2$, if $p_i$ is a
singular point of $\bar{\mathcal L}$. For $u$ sufficiently close
to $u_0$ the surface $Y_u\subset V$ and hence for such $u$ the
restriction of  \text{pr} to $Y_u$ satisfies the conditions of
Proposition \ref{dis}. It follows from flatness of the family of
branch curves $B_u$ that for almost all $u\in U$ the curves $B_u$
are reduced. According to  \cite{C-F} the restriction of a generic
projection to $Y_{u_1}$ (for fixed $u_1$) is a generic covering of
the plane and, in particular, the curve $B_{u_1}$ is cuspidal.
Therefore, if the projection is generic simultaneously for
$\mathcal P$ and for $Y_{u_1}$, then it follows from flatness of
the family $B_u$ that for almost all $u$ the curves $B_u$ have
singularities not worse than the singularities of $B_{u_1}$, that
is, for almost all $u\in U$ the curves $B_u$ are also cuspidal and
the restriction of the projection $\text{pr}$ to $Y_u$ is a
generic covering of the plane.  \qed

\subsection{Numerical data for description of a projection of a curve
$\mathcal D \subset \mathcal L$ onto the plane} Let ${\rm
pr}:\mathbb P^3\to \mathbb P^2$ be a generic projection of an
arrangement of $m$ planes $\mathcal P \subset \mathbb P^3$. Then
the curve $\bar {\mathcal L} ={\rm pr}(\mathcal L) \subset \mathbb
P^2$ has $\tau$ triple points and
\begin{equation} \label{nu}
\nu =\frac{1}{2} {m\choose 2} {m-2\choose 2}
=\frac{m(m-1)(m-2)(m-3)}{8}
\end{equation}
double points  --- points of intersection of lines ${\rm
pr}(L_{i,j})$ and ${\rm pr}(L_{k,l})$, which are projections of
skew lines $L_{i,j}$ and $L_{k,l}$ of the line arrangement
$\mathcal L \subset \mathbb P^3$.

Let $\mathcal L = \mathcal D \cup \mathcal R$ be a partition of
${m\choose 2}$ lines of $\mathcal L$ into two parts, containing
$\bar d$ and $d$ lines respectively. Denote $\bar {\mathcal D}
={\rm pr}(\mathcal D)$ and $\bar {\mathcal R} ={\rm pr}(\mathcal
R)$. Then $\nu$ double points of the curve $\bar{\mathcal L}$ are
decomposed into 3 sets: $\nu_{2}$ points, which are double on
$\bar {\mathcal D}$; $\nu_{1}$  points, which are simple on $\bar
{\mathcal D}$ (they are intersection points of $\bar {\mathcal D}$
and $\bar {\mathcal R}$); $\nu_{0}$ points, which do not lie on
$\bar {\mathcal D}$ (they are double points of $\bar {\mathcal
R}$),
\begin{equation} \label{nu1}
\nu =\nu_{2} + \nu_{1}+\nu_{0} .
\end{equation}

The curve $\bar {\mathcal D}$, which is an arrangement of $\bar d$
lines, has $\tau_{3}$ triple points and $\tau_{2} +\nu_{2}$ double
points, $\tau_{2}$ of which have being on $\mathcal D$, and
$\nu_{2}$ appeared under projection of skew lines of the curve
$\mathcal D$. Applying formula (\ref{pa}) to the curve
$\bar{\mathcal D}$, we get
$$
\frac{1}{2} (\bar d -1)(\bar d -2) =(\tau_{2} +\nu_{2}) +3\tau_{3}
-\bar d +1
$$
or $\frac{1}{2} {\bar d} (\bar d -1) =\nu_{2} +\tau_{2}
+3\tau_{3}$, from where
\begin{equation} \label{nu2}
\nu_{2} = \frac{1}{2} \bar d (\bar d -1) -\tau_{2} -3\tau_{3}.
\end{equation}

The analogous formula for the curve $\bar{\mathcal R}$ of degree
$d$ with $\tau_{0}$ triple points and $\tau_{1} +\nu_{0}$ double
points gives
\begin{equation} \label{nu0}
\nu_{0} = \frac{1}{2} d (d -1) -\tau_{1} -3\tau_{0}.
\end{equation}

\subsection{Expression of numerical characteristics of a covering
in terms of degeneracy} We expressed invariants $K_{X}^{2}$ and
$e(X)$ of an irreducible surface $X$ in terms of numerical
characteristics of a degeneracy of the surface $Y$. On the other
hand, $K_{X}^{2}$ and $e(X)$ can be expressed in terms of
numerical characteristics of a generic covering (see (\ref{k2})
and (\ref{e})). This gives an expression of numerical
characteristics of the covering in terms of numerical
characteristics of degeneration.

\begin{prop} \label{prop4}
If $Y_{u}$, $u\in U$, is a complete degeneration of an irreducible
surface $Y$, then numerical characteristics of a generic
projection $p : Y \rightarrow {\mathbb{P}}^{2}$ are expressed in
terms of numerical characteristics of the degeneration by
formulae:
\begin{equation}
 \deg B  = 2d, \label{d}
\end{equation}
\begin{equation}  g-1  = 6\tau_0+\tau_1-d, \label{g}
\end{equation}
\begin{equation} \displaystyle c  =  6\tau_0+3\tau_1,  \label{c}
\end{equation}
\begin{equation} \displaystyle n  =  4\nu_{0} .  \label{n}
\end{equation}
\end{prop}
\proof First we prove (\ref{d}), that is, that $d$ in the notation
of degree $ \deg B $ in section \ref{1.4} equals to
$d=\deg{\mathcal R}$ in the formula (\ref{L}). Let $L$ be a
generic line in ${\mathbb{P}}^{2}$, and $\bar L =p^{-1}(L)$ be the
corresponding hyperplane section of $Y$. Then $\bar L$ is an
irreducible plane curve of degree $m$ with $\bar d =\deg D$ nodes.
Hence the geometric genus $g(\bar L) =\frac{(m-1)(m-2)}{2} -\bar
d$. On the other hand,  $p: \bar L \rightarrow L$ is a covering of
degree $m$ ramified at $\deg B = (B,L)$ points. Therefore, by
Hurwitz formula, we have: $2 g(\bar L) -2 =-2m + \deg B$. It
follows that
$$
\deg B =(m-1)(m-2) -2\bar d -2 +2m =m(m-1) -2\bar d =2d .
$$

Let us prove formula (\ref{g}). From (\ref{k2}) we have $g-1
=K_{X}^{2} -9m +9d$ and by virtue of (\ref{K2})
$$
g-1 =m(m-4)^{2} +10\bar d -5\tau_{1} -6\tau_{2} -6\tau_{3} -9m +9d
.
$$
Substituting  $\bar d =\frac{(m-1)(m-2)}{2} -d$ from (\ref{L}) and
$\tau_{2} +\tau_{3}=\tau -\tau_{1} -\tau_{0}$ from (\ref{tau}), we
obtain formula (\ref{g}): $g-1 =m(m-4)^{2} +5m(m-1) -10 d -
5\tau_{1} - m(m-1)(m-2) +6\tau_{1} +6\tau_{0} -9m +9d
=6\tau_0+\tau_1-d$.

To calculate $c$ we use formula (\ref{e}): $c=-e(X)+3m+2(g-1)$.
Substituting the expression for $e(X)$ from (\ref{Etau}) and $g-1$
from (\ref{g}), we obtain:
$$c =-m^{2}(m-4) -6m -2\bar d +7\tau_{1}
+6\tau_{2}+6\tau_{3} +3m + 12\tau_{0}+2\tau_{1} -2d= $$
$$
-m^{2}(m-4) -3m -2(\bar d +d)
+6(\tau_{0}+\tau_{1}+\tau_{2}+\tau_{3}) +6\tau_{0}+3\tau_{1} .
$$
Applying formulae (\ref{L}) and (\ref{tau}), we obtain (\ref{c}).

To calculate $n$ we use the formula for genus of a plane curve
$$
\displaystyle \frac{\deg B(\deg B-3)}{2}=  g-1+c+n .
$$
From here and from formulae  (\ref{d}), (\ref{g}), and (\ref{c})
we obtain
$$
d(2d-3)=6\tau_0+\tau_1-d+ 6\tau_0+3\tau_1+n
=12\tau_0+4\tau_1+n-d,$$ that is, $n-(d-1)-12\tau_0-4\tau_1=
4\nu_{0}$ by virtue of formula (\ref{nu0}). \qed

\subsection{Degenerations of cubic surfaces. The geometric meaning of formulae
for $c$ and $n$} \label{cub} Let $\text{pr}: Y\to\mathbb P^2$ be a
generic projection of a surface $Y\subset \mathbb P^3$, $\deg
Y=3$.

It is known  \cite{G-H} that the irreducible surfaces with
ordinary singularities of degree $m=3$ are either smooth cubics or
cubics, the double curve of which is a line.

If $Y$ is a smooth cubic, then, as is known, the branch curve
$B\subset \mathbb P^2$ is a curve of degree 6 with six cusps
(which in addition lie on a conic).

If the double curve $D$ is a line ($\bar d =1$), then the surface
$Y$ has two pinches, $\omega =2$, and the branch curve $B$ is a
rational quartic ($2d =4$) with three cusps ($c=3$).

If a surface $Y$ is reducible, then either $Y=P\cup Q$ is the
union of a plane and a quadric, or $Y$ is the union of three
planes.

Now consider a  complete degeneration of a cubic $Y$ into an
arrangement of three planes $\mathcal P= P_{1}\cup P_{2}\cup
P_{3}$ with the double curve $\mathcal L= L_{1,2}\cup L_{2,3}\cup
L_{3,1}$. The arrangement $\mathcal P$ has one triple point $s=
L_{1,2}\cap L_{2,3}\cap L_{3,1}$. Let $\mathcal D$ and $\mathcal
R$ be the limit double curve and the limit ramification curve
respectively. The surface $Y$ is obtained from $\mathcal P$ by
smoothing $\mathcal R$ (smoothing outside of $\mathcal D$). We
have 4 possibilities:

1) if the generic fibre $Y$ is a smooth cubic, then $\mathcal L
=\mathcal R$, $\mathcal D= \emptyset$, that is, the double curve
$\mathcal L$ is smoothed completely; in this case from the triple
point $\bar s = p(\mathcal{T}_{0})\in \bar{\mathcal R}$ there
appears 6 cusps (and no nodes) on the branch curve $B$;

2) if the generic fibre $Y$ is an irreducible cubic,whose double
curve $D$ is a line, then $\mathcal D =L_{1,2}$ (for example), and
$\mathcal R =L_{2,3}\cup L_{3,1}$, i.e two lines in $\mathcal L$
are smoothed and one line remains double; the point
$s\in\mathcal{T}_{1}$; in this case from the point $s$ there
appear two pinches on $Y$, and the point $\bar s$ gives three
cusps on the curve $B$;

3) if the generic fibre $Y=P\cup Q$, then $\mathcal R =L_{1,2}$,
$\mathcal D =L_{2,3}\cup L_{3,1}$; the point
$s\in\mathcal{T}_{2}$; the curves $2\mathcal R$ and $\mathcal D$
are smoothed into plane conics, one of which becomes a
ramification curve $R$, and the other becomes a double curve $D$;
the curve $B$ is a conic in $\mathbb P^2$.

4) Finally, if $Y$  is a union of three planes, then $\mathcal L
=\mathcal D$, $\mathcal R= \emptyset$, and $s\in\mathcal{T}_{3}$;
the double curve $\mathcal L$ is not smoothed at all.

Now let  $\mathcal P$ be a  complete degeneration of a surface $Y$
of degree $m$ with ordinary singularities, $\mathcal D$ and
$\mathcal R$ be the limit double curve and the limit ramification
curve.  If $s \in \mathcal{T}$ a triple point on $\mathcal P$,
then (locally) in a neighborhood of the point $s$ the smoothing
outside of $\mathcal D$, or the regeneration, looks like as in the
case of the regeneration of a cubic surface, and this explains
obtained previously formulae  in the following way.

Formula (\ref{c}) is explained by the fact that the regeneration
in the case $s \in \mathcal{T}_{0}$ gives 6 cusps on the curve
$B$, and in the case $s \in \mathcal{T}_{1}$ there appears 3
cusps.

Formula (\ref{wtau}) is explained by the fact that pinches (by
two) appear only from points $s \in \mathcal{T}_{1}$.

Finally, formula  (\ref{n}) is explained as follows. Let $q$ be
one of $\nu_0$  double points of the curve $\bar{\mathcal R}$,
$q=\bar L_{i,j}\cap \bar L_{k,l}$, where $L_{i,j}$ and $L_{k,l}$
is a pair of skew lines.  In a neighborhood of a line $L_{i,j}$
the surface $\mathcal P$ is given (locally) by the equation
$z^2-x^2=0$ and smoothed surface $Y_{u}$ is given by an equation
$z^2-x^2=\varepsilon$. Projecting to the plane ${\mathbb P}^{2}$,
we obtain two branches of the branch curve $B_u$ close
("parallel") to the curve $\bar L_{i,j}$. Analogously, for the
curve $L_{k,l}$  we obtain two branches of the branch curve $B_u$
close to the line $\bar L_{k,l}$. These two pairs of branches meet
in four points. Thus, each of $\nu_{0}$ pairs of skew lines of the
curve $\mathcal R$ under regeneration gives 4 nodes on the curve
$B$, and we obtain the formula $n=4\nu_{0}$.

\section{Complete degeneracy of quartic surfaces}
\label{deg4} Let us analyze complete degenerations of surfaces of
degree $m=4$. If a surface $Y$ is reducible, then it is obvious
that complete degenerations of components of $Y$ give a complete
degeneration of $Y$. Thus, we can assume the surface $Y$ to be
irreducible since the case $m=2$ is trivial, and the case $m=3$
has been described in section \ref{cub}.

\subsection{Degenerations of irreducible quartics with
ordinary singularities} \label{4.1} The description of all
irreducible quartics $Y\subset\mathbb{P}^{3}$ with ordinary
singularities can be found in \cite{G-H}. There are 6 types of
such surfaces in accordance with the types of the double curve
$D\subset Y$:

1)  $Y$  is smooth, that is, $D=\emptyset$;

2) $D$ is a line;

3) $D$ is a plane conic;

4) $D$ is a pair of skewed lines;

5) $D$ is the union of three lines meeting at a point;

6) $D$ is a rational normal curve of degree $3$.

We show that in each of these cases there is a complete
degeneration of $Y$. Let $F(x)=0$ be an equation of $Y$, and
$H_{i}(x)=0$, $i=1,\dots ,4$, be equations of four planes
$P_{1},\dots , P_{4}$ in general position in $\mathbb{P}^{3}$.

In cases 1), 2), 4) and 5) we can obtain the desired complete
degenerations  $Y_{u}$, $u\in {\mathbb{C}}$, in a form
$$
u\cdot F(x) + (1-u)H_{1}(x)H_{2}(x)H_{3}(x)H_{4}(x) =0 ,
$$
where in case 1) linear functions $H_{1}(x), \dots , H_{4}(x)$ are
arbitrary and such that $P_{1},\dots , P_{4}$ are in general
position; in case 2) $H_{1}(x)=H_{2}(x)=0$ are equations of the
line $D$; in case 4) $H_{1}(x)=H_{2}(x)=0$ and
$H_{3}(x)=H_{4}(x)=0$ are equations of skew lines forming $D$; in
case 5) $H_{1}(x)H_{2}(x)H_{3}(x) =0$ is an equation of the union
of three planes, the double curve of which is the union of three
lines forming $D$. In all cases the double curve $\mathcal L$ of
the arrangement $\mathcal P$ of four planes $P_{1},\dots , P_{4}$
contains the double curve $D_{u}\subset Y_{u}$ and the limit
double curve $\mathcal D$ coincides with $D$.

In the case 3) the double curve $D$ is a complete intersection of
a plane and a quadric. We consider complete degenerations with a
complete intersection double curve in general in section
\ref{cmlt}.

Consider case 6). As is known, all smooth space  curves of degree
$3$ are projectively equivalent. Such a curve $D$ has
parametrization $x=t$, $y=t^{2}$, $z=t^{3}$ in appropriate affine
coordinates. It is not a complete intersection and is defined by
three equations: $y=x^{2}, xy=z, y^{2}=xz$. In homogeneous
coordinates $(x_{1}:x_{2}:x_{3}:x_{4}) =(x:y:z:1)$ in
$\mathbb{P}^{3}$ the curve $D$ is the intersection of three
quadrics
$$
x_{2}x_{4}= x_{1}^{2}, \;  x_{1}x_{2}= x_{3}x_{4}, \;
x_{2}^{2}=x_{1}x_{3} .
$$

Consider a family of curves $D_{u}$: $x=t$, $y=ut^{2}$,
$z=u^{3}t^{3}$. The curves $D_{u}$ are given by three equations
$y=ux^{2}, u^{2}xy=z, uy^{2}=xz$, or in homogeneous coordinates
$$
x_{2}x_{4} -u x_{1}^{2} =0, \; x_{3}x_{4} - u^{2} x_{1}x_{2}= 0 ,
 \; x_{1}x_{3} -u x_{2}^{2}= 0.
$$
The family $D_{u}$ defines a degeneration of the curve $D$ (for
$u=1$) to the curve $\mathcal D$ (for $u=0$), which is given by
equations: $x_{2}x_{4}=0 ,\, x_{3}x_{4}=0 ,\, x_{1}x_{3}=0$. If
$P_{i}$ is a plane given by equation $x_{i}=0$, $L_{i,j}=P_{i}\cap
P_{j}$, then $\mathcal D$ is a chain of lines $L_{1,4}\cup
L_{4,3}\cup L_{3,2}$.

Consider a family of surfaces $Y_{u}$ given by equation
$$
u(x_{2}x_{4} -u x_{1}^{2})(x_{3}x_{4} - u^{2} x_{1}x_{2}) +
u(x_{3}x_{4} - u^{2} x_{1}x_{2})(x_{1}x_{3} -u x_{2}^{2}) +
$$
$$
(x_{2}x_{4} -u x_{1}^{2})(x_{1}x_{3} -u x_{2}^{2}) = 0 .
$$
For $u\neq 0$ the double curves are smooth space cubic curves
$D_{u}$, and for $u=0$ the surface $Y$ degenerates to a plane
arrangement $Y_{0}=\{ x_{1}x_{2}x_{3}x_{4} =0\}$ and the limit
double curve is $\mathcal D$.

\subsection{Description of possible arrangements of limit double curves}
In case  $m=4$ the graph $\Gamma (\mathcal L)$ consists of the
union of all edges of a tetrahedron. By Lemma \ref{R}, a curve
$\mathcal D\subset \mathcal L$ can be a limit double curve of a
degeneration of an irreducible quartic with ordinary singularities
if the graph $\bar \Gamma (\mathcal R)$ is connected. The graph
$\bar \Gamma (\mathcal R)$ is obtained from $\Gamma (\mathcal L)$
by removing some edges corresponding to lines $L_{i,j}\subset
\mathcal D$. It is obvious that the graph $\bar \Gamma (\mathcal
R)$ remains connected only in the following cases: 1) we remove
nothing from  $\Gamma (\mathcal L)$; 2) we remove one edge from
$\Gamma (\mathcal L)$; 3) we remove two adjacent edges; 4) we
remove two skew edges; 5) we remove three edges going out from one
vertex; 6) we remove a chain of three edges. Removing of four
edges leads to a non isolated vertex, that is, to a non connected
graph $\bar \Gamma (\mathcal R)$.

Thus, up to renumbering of planes of $P_{1},\dots , P_{4}$, the
curve $\mathcal D$ is one of the following:

1) $\mathcal D=\emptyset$;

2) $\mathcal D=L_{1,2}$;

3) $\mathcal D=L_{1,2}\cup L_{1,3}$;

4) $\mathcal D=L_{1,2}\cup L_{3,4}$;

5) $\mathcal D=L_{1,2}\cup L_{1,3}\cup L_{2,3}$;

6) $\mathcal D=L_{1,2}\cup L_{2,3}\cup L_{3,4}$.

\noindent In section  \ref{4.1} we showed that all such $\mathcal
D$ are realized as limit double curves. Thus, we obtain the
following
\begin{thm} For any surface of degree $4$ with ordinary
singularities in $\mathbb P^3$ there exists a complete
degeneration.

If $\mathcal P$ is an arrangement of four planes in general
position in $\mathbb P^3$, $\mathcal D \subset \mathcal L$ is any
line arrangement, then there exists a complete degeneration of
surfaces of degree $4$ with ordinary singularities, for which
$\mathcal D$ is the limit double curve.
\end{thm}

\section{Existence of line arrangements which are not limit}
\label{m5} In the previous section we showed that in the case of
surfaces with ordinary singularities of degree $m\leq 4$ the matters
go as well as in the case of plane curves. Responses to two both
questions, put in the introduction, are affirmative.  But if $m\geq
5$, the matters do not stand so good. In this section we show that
there are at least seven arrangements of double lines $\mathcal D$
in plane  arrangements $\mathcal P$, $\deg \mathcal P=5$, such that
$\mathcal P$ cannot be smoothed outside of of $\mathcal D$.

\begin{figure}[h]
\begin{picture}(0,160)(160,-130)

\put(30,18){\line(0,-1){30}} \put(30,18){\circle*{3}}
\put(30,-12){\circle*{3}} \put(60,18){\line(0,-1){30}}
\put(60,18){\circle*{3}} \put(60,-12){\circle*{3}}
\put(30,18){\line(1,0){30}} \put(30,-12){\line(1,0){30}}
\put(30,18){\line(1,-1){30}} \put(60,18){\line(-1,-1){13}}
\put(30,-12){\line(1,1){13}} \put(30,-30){$\Gamma_{0,0,4}^{6,4}$}

\put(100,18){\line(0,-1){30}} \put(100,18){\circle*{3}}
\put(150,18){\circle*{3}} \put(100,-12){\circle*{3}}
\put(130,18){\line(1,0){20}} \put(130,18){\circle*{3}}
\put(130,-12){\circle*{3}} \put(100,18){\line(1,0){30}}
\put(100,-12){\line(1,0){30}} \put(100,18){\line(1,-1){30}}
\put(130,18){\line(-1,-1){13}} \put(100,-12){\line(1,1){13}}
\put(110,-30){$\Gamma_{1,1,3}^{6,2}$}

\put(185,13){\line(0,-1){20}} \put(185,13){\circle*{3}}
\put(185,-7){\circle*{3}} \put(195,3){\circle*{3}}
\put(185,13){\line(1,-1){10}} \put(185,-7){\line(1,1){10}}
\put(195,3){\line(1,0){15}} \put(210,3){\circle*{3}}
\put(210,3){\line(1,0){15}} \put(225,3){\circle*{3}}
\put(190,-30){$\Gamma_{1,3,1}^{5,1}$}

\put(255,13){\circle*{3}} \put(255,-7){\circle*{3}}
\put(265,3){\circle*{3}} \put(255,13){\line(1,-1){10}}
\put(255,-7){\line(1,1){10}} \put(265,3){\line(1,0){15}}
\put(280,3){\circle*{3}} \put(280,3){\line(1,0){15}}
\put(295,3){\circle*{3}} \put(260,-30){$\Gamma_{3,1,1}^{4,0}$}

\put(30,-47){\line(0,-1){20}} \put(30,-47){\circle*{3}}
\put(30,-67){\circle*{3}} \put(40,-57){\circle*{3}}
\put(30,-47){\line(1,-1){10}} \put(30,-67){\line(1,1){10}}
\put(40,-57){\line(1,0){15}} \put(55,-57){\circle*{3}}
\put(32,-90){$\Gamma_{1,2,1}^{4,1}$}

\put(150,-47){\circle*{3}} \put(150,-67){\circle*{3}}
\put(160,-57){\circle*{3}} \put(150,-47){\line(1,-1){10}}
\put(150,-67){\line(1,1){10}} \put(160,-56){\line(1,0){15}}
\put(175,-57){\circle*{3}} \put(150,-90){$\Gamma_{3,0,1}^{3,0}$}

\put(255,-57){\circle*{3}} \put(255,-57){\line(1,0){15}}
\put(270,-57){\circle*{3}} \put(270,-57){\line(1,0){15}}
\put(285,-57){\circle*{3}} \put(265,-67){\circle*{3}}
\put(265,-67){\line(1,0){15}} \put(280,-67){\circle*{3}}
\put(259,-90){$\Gamma_{(2)(2,1)}^{(1,0)(2,0)}$}

\put(150,-120){$\text{\rm Fig}.\, 2$}
\end{picture}
\end{figure}

\begin{thm} \label{thm-5}
The line arrangements $\mathcal D $ with graphs, depicted on Fig.
{\rm 2}, are not limit for complete degenerations of surfaces with
ordinary singularities of degree {\rm 5}.
\end{thm}
\proof Assume that line arrangements with graphs, depicted on Fig.
{\rm 2}, are  limit double curves for complete degenerations of
surfaces $Y$ of degree $m =5$ with ordinary singularities.

In all cases the graph $\bar\Gamma(\mathcal R)$ is connected and
by Lemma \ref{R} the surface $Y$, for which $\mathcal D$ is the
limit double curve, is irreducible.

We can obtain numerical characteristics of curves $\mathcal D$. In
our case $m=5$ and therefore, $d+\bar d ={5\choose 2}$, $\tau
={5\choose 3}$. The number $\tau_{3}$ is equal to the number of
triangles in $\Gamma(\mathcal D)$. To calculate $\tau_{2}$ we use
the formula (\ref{v}); $\tau_{1}$ and $\tau_{0}$ are obtained from
formulae (\ref{tau123}), (\ref{tau210}), and (\ref{tau}):
$$
\tau_{1} +2\tau_{2} +3\tau_{3} =3\bar d , \;  \tau_{2} +2\tau_{1}
+3\tau_{0} =3d , \; \tau_{0} +\tau_{1} +\tau_{2}+\tau_{3}=\tau.
$$
By (\ref{g}), (\ref{n}), and (\ref{c}) we can find numerical
characteristics of the double curve $D\subset Y$ and of the branch
curve $B$: the numbers $\bar g$, $t$, $k$ we find from Proposition
\ref{prop1}, and the numbers $g$, $c$, $n$ we find from
Proposition \ref{prop4} applying formulae (\ref{g}), (\ref{c}),
(\ref{n}), and (\ref{nu0}).

Consider each of the line arrangements $\mathcal D $ with graphs,
depicted in Fig. {\rm 2}, and show that in all cases we obtain a
contradiction with assumption made above for the arrangements to be
limit.

For a line arrangement of type $\Gamma_{3,0,1}^{3,0}$ we have:
$\bar d =3$, $\tau_{3}=0$ and $\tau_{2}=3$. By Proposition
\ref{prop1}, the double curve $D\subset Y$ is an irreducible
smooth curve of degree $\bar d =3$ and genus $\bar g =1$. Hence,
$D$ is a plane cubic. Let $D$ lies in a plane $P$. Then for any
line $L\subset P$ the intersection number $(L, Y)_{\mathbb P^2}
\geq 6$, and since $\deg Y =5$, we have $L\subset Y$. Thus,
$P\subset Y$ and the surface $Y$ is reducible. A contradiction.
The same arguments can be applied for line arrangements of the
following two types.

For a line arrangement of type $\Gamma_{(2)(2,1)}^{(1,0)(2,0)}$ we
have: $\bar d =3$, $\tau_{3}=0$ and  $\tau_{2}=1$. By Proposition
\ref{prop1}, we have $k=2$, $\bar g =0$ and the double curve
$D\subset Y$ is a disconnected union of a line $L_{1}$ and a plane
conic $Q$. Let $Q$ lies in a plane $P$ and the point $A=L_{1}\cap
P$. As in the previous case, we see that any line $L\subset P$,
passing through $A$, lies in the surface $Y$ and, consequently,
$P\subset Y$. Again we obtain a contradiction with the
irreducibility of $Y$.

For a line arrangement of type $\Gamma_{1,2,1}^{4,1}$ we have:
$\bar d =4$, $\tau_{3}=1$ and $\tau_{2}=2$. By Proposition
\ref{prop1}, we have $k=2$, $\bar g =0$, $t=1$. Hence,  $D=L\cup
C$ consists of two irreducible components, where $L$ is a line,
and $C$ is a curve of degree 3 with one node. Consequently, the
curve $C$ lies in a plane $P$. As in the previous cases we obtain
that $P\subset Y$ and this gives a contradiction with the
irreducibility of $Y$.

We consider a line arrangement of type $\Gamma_{3,1,1}^{4,0}$ in
Theorem \ref{nonreal}, where it is shown that the curve $\mathcal
D$ can not be the limit double curve of a complete degeneration of
a surface $Y$ of any degree $m$, and not only $m=5$.

For a line arrangement of type $\Gamma_{1,3,1}^{5,1}$ we have:
$\bar d =5$, $\tau_{3}=1$, $\tau_{2}=3$, and $\tau_{1}=6$. By
(\ref{K2}) and (\ref{E}), we obtain $K^2_X=1$ and $e(X)=-1$. It is
known from the classification of algebraic surfaces that if
$K^2_X>0$, then $e(X)>0$. Thus such surfaces do not exist.

A line arrangement of type $\Gamma_{0,0,4}^{6,4}$ is a special
case (for $m=5$) of line arrangements which are considered below
in Proposition \ref{nonreal2}.

For a line arrangement of type $\Gamma_{1,1,3}^{6,2}$ we have:
$\bar d =6$, $d=4$, $\tau_{3}=2$, $\tau_{2}=4$, $\tau_{1}=4$,
$\tau_{0}=0$. The branch curve $B$ is cuspidal and have numerical
characteristics: $\deg B =8$, $g=1$, $c=12$, $n=8$. Let us show
that the curve $\hat{B}$ dual to $B$ also is cuspidal.  By
Pl\"{u}cker formula, we obtain $ \deg \hat{B} =4$. Since
$g(\hat{B})=1$ and a curve of degree 4 has arithmetic genus
$p_{a}(\hat{B}) =3$, then  $\hat{B}$ has either two singular
points with $\delta =1$, or one singular point with $\delta =2$.
The Milnor number $\mu$ and $\delta$ are connected by formula $\mu
= 2\delta -r+1$, where $r$ is the number of branches.
Singularities with $\delta =1$ are either singularities of type
$A_1$ -- nodes, or of type $A_2$ -- cusps, therefore, in the first
case the curve $\hat{B}$ is cuspidal. The second case is
impossible, since the singularity dual to $A_{3}$ is $A_{3}$, and
dual to $A_{4}$ is $E_{8}: x^{3}+y^{5}=0$. But a curve dual to
$\hat{B}$ is the curve $B$, which has only nodes and cusps. Thus,
the curve $\hat{B}$ is cuspidal.  In this case it follows from
Pl\"{u}cker formula  that $3\deg B-c=3\deg \hat{B}-\hat{c}$, where
$\hat{c}$ is the number of cusps of $\hat{B}$. Therefore,
$\hat{c}=0$ and, consequently, the curve $\hat{B}$ is nodal and
the number of node of $\hat{B}$ is equal to $\hat{n}=2$. For such
a curve $\hat{B}$ there exists a generic covering of the plane of
degree $4=\deg \hat{B}$, branched along $B$ (see \cite{Ku}). On
the other hand, a generic projection of a surface with ordinary
singularities  in $\mathbb P^{3}$ onto $\mathbb P^2$ is a generic
covering, in our case of degree $5$. As is shown in \cite{Ku3},
generic coverings, which are generic projections, are uniquely
defined by their branch curves always, except the case of surfaces
of degree $4$ with singularities consisting of three lines
intersecting in a point.\footnote{In \cite{Ku3} this statement was
formulated only for surfaces obtained as a generic linear
projection into $\mathbb P^3$ of smooth surfaces imbedded into
some $\mathbb P^N$.  But the proof of this statement used only the
assumption that surfaces in $\mathbb P^3$ have only ordinary
singularities.} We obtain a contradiction with the uniqueness of
the covering.  \qed

\section{Existence of surfaces not possessing complete degenerations}
\label{examp}

In this section we show that there exist surfaces in $\mathbb P^3$
with ordinary singularities which can not be degenerated into
plane arrangements in general position.

\subsection{Zeuthen's problem} \label{example1} The problem of existence of a
complete degeneration for any surface with ordinary singularities
is closely connected with the problem of degeneration of its
double curve $D$ into a line arrangement.  In the case of a smooth
space curve $D$ it is a famous Zeuthen's problem. In \cite{Ha}, a
negative solution of Zeuthen's problem is obtained, that is, it is
shown that there exist smooth projective curves $D\subset \mathbb
P^3$, which can not be degenerated into a line arrangement with
double points.  One of such curves has degree $30$ and genus
$113$. We call it {\it Hartshorne's curve}.

\begin{thm} \label{Hart} There exist surfaces in $Y\subset \mathbb P^3$
with ordinary singularities, which can not be degenerated into a
plane arrangement in general position.
\end{thm}

\proof Let $D$ be the Hartshorne's curve. In the next subsection
we prove (Proposition \ref{exist}) that there exists a surface
$Y\subset \mathbb P^3$ with ordinary singularities, the double
curve of which is $D$. By \cite{Ha}, the curve $D$  can not be
degenerated into a line arrangement with double points and,
consequently, $Y$ can not be degenerated into a plane arrangement
in general position.  \qed

\subsection{Existence of a surface with ordinary singularities, the double
curve of which is any smooth curve} Let $D$ be a smooth (nor
necessary irreducible) curve in $\mathbb P^3$. Then for any point
$x\in D$ there exist a Zariski open set $U_x$ in $\mathbb P^3$ and
two rational functions $f_{x,1}$ and $f_{x,2}$, regular in $U_x$
and such that $f_{x,1}$ and $f_{x,2}$ are local parameters in
$U_x$ and the curve $D\cap U_x$ is defined by equations
$f_{x,1}=f_{x,2}=0$. Let $f_{x,j}=\frac{F_{x,j}}{G_{x,j}}$, where
$F_{x,j}$ and $G_{x,j}$ are relatively prime homogeneous
polynomials of degree $M_{x,j}$. The curve $D$ being compact, we
can choose a finite covering of $D$ by open sets $U_{x_1},\dots,
U_{x_k}$ such that polynomials $F_{x_i,j}$, $1\leq i\leq k$,
$j=1,2$, generate the homogeneous ideal of $D$. Set
\begin{equation} \label{M}
M(D)=2\max_{1\leq i\leq k}(\max_{j=1,2}M_{x_i,j})+1.\end{equation}

\begin{prop} \label{exist} For any smooth projective curve
$D\subset \mathbb P^3$ and for any natural number $m\geq M=M(D)$,
where $M(D)$ is the number defined in {\rm (\ref{M})}, there
exists a projective surface $Y$ with ordinary singularities of
degree $m$ in $\mathbb P^3$ such that $D=\text{Sing} \, Y$.
\end{prop}

\proof Let $\sigma: \widetilde{\mathbb P}^3\to \mathbb P^3$ be a
monoidal transformation with center at $D$, $E=\sigma^{-1}(D)$ its
exceptional divisor and $\widetilde P=\sigma^*(P)$ the total
preimage of the plane $P\subset \mathbb P^3$.

The surfaces $Y_{i,j}=\{ F_{x_i,j}=0\}\subset \mathbb P^3$, $\deg
Y_{i,j}=\deg F_{x_i,j}=M_{x_i,j}$, $j=1,2$, meet transversally
along $D$ in $U_{x_i}$. It is obvious that for any plane $P$ in
$\mathbb P^3$ the divisors
$$Y_{i,1}+Y_{i,2}+(m- M_{x_i,1}-M_{x_i,2})\widetilde P$$
are zeroes of sections of the sheaf $\mathcal
O_{\widetilde{\mathbb P}^3}(m\widetilde P-2E)$ for $m\geq M$. From
this it is easy to see that the linear system $\mid m\sigma^*
(P)-2E\mid$ is not empty, does not have fixed components and base
points, and for any points $x ,y\in \widetilde{\mathbb P}^3$ there
exists a divisor of this system such that it does not go through
the point $y$ and is a smooth reduced surface at $x$. According to
Bertini theorem the generic member $X$ of the linear system $\mid
m\sigma^* (P)-2E\mid$ is a smooth surface, and it is easy to see
that its image $Y=\sigma(X)$ is a surface in $\mathbb P^3$ of
degree $m$ with ordinary singularities along $D$. \qed

\section{Complete degeneracy of surfaces with complete intersection double curves}
\label{cmlt}

As has been shown in the previous section, in general the answer
to the question about existence of a complete degeneration is in
negative even in the case of a smooth double curve. In this
section we show that the answer to this question is affirmative in
the case when the double curve is complete intersection.

\subsection{Equations of surfaces with complete intersection double curves}
The following proposition gives a description of equations of
surfaces, double curves of which are complete intersections.
\begin{prop} \label{compl}
Let an irreducible surface $Y\subset \mathbb P^3$ with ordinary
singularities, $\deg Y=m$, given by equation $F(x)=0$, has a
smooth double curve $D=Y_{1}\cap Y_{2}$, which is a complete
intersection of surfaces $Y_{1}$ and $Y_{2}$ of degrees $m_{1}$
and $m_{2}$. Then the polynomial $F$ can be written in a form
\begin{equation} \label{F}
F=AF_{1}^{2}+BF_{1}F_{2}+CF_{2}^{2} ,
\end{equation}
where $F_{1}(x)=0$ and $F_{2}(x)=0$ are equations of surfaces
$Y_{1}$ and $Y_{2}$, and $A$, $B$, and $C$ are homogeneous
polynomial.

Conversely, if $F$ is written in the form {\rm (\ref{F})} and
$m\geq 2m_{1}+1$ {\rm (}let $m_{1}\geq m_{2}${\rm )}, and
polynomials $A$, $B$, and $C$ are sufficiently general, then the
surface $Y$ has only ordinary singularities and its double curve
is a complete intersection.
\end{prop}
\proof Since the curve $D$ is a complete intersection, the
homogeneous ideal $I(D)$ is generated by two elements, $I(D)
=(F_{1},F_{2})$. Write $F$ in a form $F=K_{1}F_{1}+ K_{1}F_{1}$.
It follows from transversality of intersection of surfaces $F_1=0$
and $F_2=0$ along $D$ that differentials $dF_1$ and $dF_2$ are
linear independent at each point of $D$. In addition, the
differential $dF$ equals to zero at each point of $D$, since $D$
is the double curve of $Y$. At each point of $D$ we have
$dF=K_1dF_1 +K_2dF_2$ and it follows from linear independence of
differentials $dF_1$ and $dF_2$ at these points that $K_1$ and
$K_2$ belong to the ideal $I(D)$, that is, $K_1=S_1F_1+S_2F_2$ and
$K_2=S_3F_1+S_4F_2$, where $S_i$ are some homogeneous polynomials.
Substituting these expressions for $K_1$ and $K_2$ to the
expression for $F$, we obtain the desired form of $F$. The second
part of lemma follows from Proposition \ref{exist} and the
estimate (\ref{M}). \qed

\subsection{Construction of a degeneration}
We prove the following theorem.
\begin{thm} \label{com-int} Surfaces with ordinary
singularities, double curves of which are complete intersections,
have complete degenerations.
\end{thm}
\proof Let the double curve $D$ of a surface $Y$ be given by
equations $F_1=0$, $F_2=0$, $\deg F_1=m_1$, $\deg F_2=m_2$. Then,
by Proposition \ref{compl}, the surface $Y$ is given by an
equation of the form
$$AF_1^2+BF_1F_2+CF_2^2=0.$$
Consider a deformation of this equation killing the first and the
third terms, and transforming the second term to a product of
linear forms. Let a family of surfaces $Y_u$, $u\in \mathbb C$, be
given by an equations
\begin{equation} \label{family}
\begin{array}{l} (uB+(1-u)\bar B)(uF_1+(1-u)\bar F_1)(uF_2+(1-u)\bar F_2)+
\\ uA(uF_1+(1-u)\bar F_1)^2+uC(F_2+(1-u)\bar F_2)^2=0,\end{array}\end{equation}
where $\bar F_1=H_{1}\cdot .\, .\, .\cdot  H_{m_1}$, $\bar
F_2=H_{m_1+1}\cdot .\, .\, .\cdot H_{m_1+m_2}$, $\bar
B=H_{m_1+m_2+1}\cdot .\, .\, .\cdot H_{m}$ are products of linear
forms such that zeroes of the forms are planes $P_1,\dots, P_m$ in
general position.

It is easy to see that for $u=1$ the surface $Y_u$ coincides with
$Y$. It is obvious also that for values $u$, close to $1$, curves
$D_u$, given by equations
$$uF_1+(1-u)\bar F_1=uF_2+(1-u)\bar F_2=0,$$ are smooth and they are
double curves of surfaces $Y_u$ with ordinary singularities. For
$u=0$ the degenerated fibre is given by an equation $H_1\cdot .\,
.\, .\cdot H_m=0$ and the limit double curve $D_0$ is given by
equations
\begin{equation} \label{double} H_{1}\cdot
.\, .\, .\cdot H_{m_1}=H_{m_1+1}\cdot .\, .\, .\cdot
H_{m_1+m_2}=0.\end{equation}

In particular, for $m_1=m_2=2$ the graph of the limit double curve
$\mathcal D=D_0$ is depicted in Fig. 3.

\begin{figure}[h]
\begin{picture}(0,60)(180,-35)
\put(170,18){\line(0,-1){30}} \put(170,18){\circle*{3}}
\put(170,-12){\circle*{3}} \put(200,18){\line(0,-1){30}}
\put(200,18){\circle*{3}} \put(200,-12){\circle*{3}}
\put(170,18){\line(1,0){30}} \put(170,-12){\line(1,0){30}}
\put(173,-40){$\text{\rm Fig}.\, 3$}
\end{picture}
\end{figure}

\section{Potentially limit and absolutely not limit line arrangements}
 \label{Pot}

Until now we have considered complete degenerations of  surfaces
with ordinary singularities for fixed degree $m$. The line
arrangements $\mathcal D \subset {\mathcal L}_{m}$  which can be
obtained in the degenerate fibre are called now $m$-{\it limit
curves}. (Here ${\mathcal L}_{m}$ is the double curve of an
arrangement of $m$ planes in general position). In section
\ref{m5} it was shown that there exist line arrangements $\mathcal
D$, which are not $5$-limit. We want to investigate the dependence
on $m$ of a curve $\mathcal D$ to be $m$-limit.

If for some $m_{0}$ a curve $\mathcal D$ is not $m_{0}$-limit ,
but for sufficiently big $m$ and for an embedding $\mathcal D
\subset {\mathcal L}_{m}$, the curve $\mathcal D$ is $m$-limit,
then $\mathcal D$ is called {\it potentially limit}. If for any
embedding $\mathcal D \subset {\mathcal L}_{m}$ a curve $\mathcal
D$ is not $m$-limit, then such a line arrangement $\mathcal D$ is
called {\it absolutely not limit}.

\subsection{Examples of potentially limit line arrangements}
\label{pot1} Let us show that a line arrangement $\mathcal D$ of
type $\Gamma_{3,0,1}^{3,0}$ is potentially limit. It was shown in
Theorem \ref{thm-5} that $\Gamma_{3,0,1}^{3,0}$ is not 5-limit. Let
us show that $\Gamma_{3,0,1}^{3,0}$ is $m$-limit for $m\geq 7$. The
line arrangement $\mathcal D$ consists of three lines $L_{1} \cup
L_{2} \cup L_{3}$, which lie in a plane $P$ and are cut out by three
planes defined by equations $H_{1}=0$, $H_{2}=0$, and $H_{3}=0$.
Such an arrangement $\mathcal D$ is limit for surfaces $Y$, double
curves $D$ of which are complete intersections of a smooth cubic and
a plane. Indeed, let $F_{1}(x)=0$ be an equation of a smooth cubic
$Y_{1}$, and $F_{2}(x)=0$ be an equation of a plane $Y_{2}=P$.
Consider a surface $F(x)=0$, where $F$ is defined by formula
(\ref{F}) and consider a family of surfaces $Y_{u}$ defined by
equation (\ref{family}). It follows from (\ref{double}) that the
limit double curve $\mathcal D$ is defined by equations
$H_{1}H_{2}H_{3}=F_2=0$ and it is of type $\Gamma_{3,0,1}^{3,0}$.

Analogously one can show that a line arrangement of type
$\Gamma_{(2)(2,1)}^{(1,0)(2,0)}$ also is potentially limit.

\subsection{An example of absolutely not limit line arrangement}
\label{pot2} The following theorem gives an example of absolutely
not limit line arrangement. Let us note that this line arrangement
is a degeneration of a smooth space curve of degree $4$ and genus
$1$.
\begin{thm} \label{nonreal} The line arrangement
$\mathcal D$ of type $\Gamma_{3,1,1}^{4,0}$ {\rm (}see Fig.
$2${\rm )} is absolutely not limit.
\end{thm}
\proof Let $\mathcal D$ be $m$-limit for some $m$, that is,
$\mathcal D$ is the limit double curve of a complete degeneration
of surfaces $Y_{u}$, $u\in U$, defined by equations $F_{u}(x)=0$.
Calculations analogous to calculations in the proof of Theorem
\ref{thm-5}, give: $\bar d =4$, $\tau_{3}=0$, and $\tau_{2}=4$. By
Proposition \ref{prop1}, the double curve $D\subset Y$ is a smooth
irreducible curve in $\mathbb{P}^{3}$ of degree $\bar d =4$ and
genus $\bar g =1$.  As is known, such curves $D$ are complete
intersections of two quadrics, and, consequently, a polynomial $F$
defining the surface $Y$ can be written in the form (\ref{F}):
\begin{equation} \label{Q}
F=AQ_{1}^{2}+BQ_{1}Q_{2}+CQ_{2}^{2} ,
\end{equation}
where $Q_{1}(x)=0$, $Q_{2}(x)=0$ are equations of these quadrics,
and $A$, $B$, and $C$ are polynomials of degree $m-4$.

Let us show that the family of surfaces $F_{u}(x)=0$, $u\in U$,
can be written in a form
\begin{equation} \label{Qu}
F_{u}=A_{u}Q_{1,u}^{2}+B_{u}Q_{1,u}Q_{2,u}+C_{u}Q_{2,u}^{2}
\end{equation}
(may be, after base change). For this let us consider the universal
family of surfaces given by equations of the form (\ref{Q}). The
base of this family $\mathcal F$ is an open subset in the space of
coefficients of forms $Q_{1}$, $Q_{2}$, $A$, $B$, $C$. Denote by
${\mathcal H}_{m,4,1}$ the space parametrizing the surfaces of
degree $m$, the double curve of which is a smooth curve of degree 4
and genus 1. Obviously, we have a rational dominant map $\mathcal F
\rightarrow {\mathcal H}_{m,4,1}$. The family of surfaces
$F_{u}(x)=0$ defines a map $U\rightarrow {\mathcal H}_{m,4,1}$. We
can assume that $U \subset {\mathcal H}_{m,4,1}$. If a curve $\tilde
U \subset \mathcal F $ is mapped to the curve $U$ we get a family of
surfaces (\ref{Qu}) parametrized by points of $\tilde U$.

A line  arrangement $\mathcal D$ of type $\Gamma_{3,1,1}^{4,0}$
consists of three lines in a plane and a fourth line not in this
plane and intersecting one of these lines.  Such curves are
degenerations of space elliptic curves of degree 4 (see
\cite{Got}). But, if a curve $\mathcal D$ is the limit double
curve of the family (\ref{Qu}), then the double curves $D_{u}$ are
given by a family of ideals $J_u=(Q_{1,u},Q_{2,u})$. By
\cite{Got}, the degenerate curve to coincide with $\mathcal D$ it
is necessary for $u=u_0$ the quadratic forms $Q_{1,u_0}$ and
$Q_{2,u_0}$ to split into a product of linear forms with one
common form: $Q_{1,u_0}=HH_1$ and $Q_{2,u_0}=HH_2$. But, then it
follows from (\ref{Qu}) that the degenerate surface $Y_{u_0}$
contains a multiple plane $H=0$, and we obtain a contradiction
with the definition of a complete degeneration of surfaces with
ordinary singularities. \qed

\section{Virtual degeneracy}
If for a surface $Y\subset \mathbb P^3$ of degree $m$ with
ordinary singularities there exists a complete degeneration in the
sense of the definition given in introduction, then, for brevity,
we call $Y$ {\it a completely degenerative surface}. In section
\ref{examp} we gave examples of surfaces, which are not completely
degenerative.

\subsection{Necessary conditions for complete degeneracy}
Let us weaken the notion of degeneracy of a surface. Recall that
in sections \ref{type1} and \ref{type2} we defined the type of an
irreducible surface $Y\subset \mathbb P^3$ of degree $m$ with
ordinary singularities and the type of a pair $(\mathcal P
,\mathcal D)$ as collections of numerical data:
$${\rm type\,}(Y) =(m,\bar d, k, \bar g, t), \qquad {\rm
type\,}(\mathcal P ,\mathcal D) =(m,\bar d,
k,\tau_{2},\tau_{3}).$$ By Proposition \ref{prop1}, if $Y$ is a
completely degenerative surface, then the types  ${\rm type\,}(Y)$
and ${\rm type\,}(\mathcal P ,\mathcal D)$ define each other. We
call collections of numbers $(m,\bar d, k, \bar g, t)$ and
$(m,\bar d, k,\tau_{2},\tau_{3})$ {\it corresponding} to each
other if $\tau_{3}=t$, and numbers $\tau_{2}$ and $\bar g$ are
connected by formula (\ref{tau2}): $\tau_{2} =  \bar d + \bar g
-k$. We call an irreducible surface $Y$ {\it virtually
degenerative} if there exists an irreducible pair $(\mathcal P
,\mathcal D)$ the type of which corresponds to the type of the
surface $Y$. Thus, a surface $Y$ can be not degenerative by a
trivial reason: there does not exist a pair $(\mathcal P ,\mathcal
D)$ of corresponding type.

Analogously, we call a pair $(\mathcal P ,\mathcal D)$ {\it
virtually smoothable outside of $\mathcal D$} (or we call a line
arrangement $\mathcal D$ {\it virtually} $m$-{\it limit}) if there
exists an irreducible surface $Y$ of degree $m$ which has a type
corresponding to the type of the pair $(\mathcal P ,\mathcal D)$.
Thus, the term "virtual" says about fulfilment of necessary
numerical conditions for existence of a complete degeneration.

The proof of Proposition \ref{prop4} does not use essentially the
fact of complete degeneracy of $Y$, but uses only arising from it
connection between the type of $Y$ and the type of the limit pair
$(\mathcal P,\mathcal D)$, given by formulae (\ref{tau3}):
$t={\tau}_{3}$,  and (\ref{tau2}): $\overline g=\tau_2- \overline
d+k$. Therefore, in the case of virtual degeneracy an analog of
Proposition \ref{prop4} holds.

\begin{prop} \label{propx}
Let the type $(m,\bar d, k, \bar g, t)$ of a surface $Y$
corresponds to the type $(m,\bar d, k,\tau_{2},\tau_{3})$ of a
pair $(\mathcal P ,\mathcal D)$. Then the same formulae {\rm
(\ref{d}) -- (\ref{n})}, as in Proposition {\rm \ref{prop4}}, hold
for numerical characteristics of a generic projection $p:Y\to
\mathbb P^2$.
\end{prop}

\subsection{Examples of virtually $m$-limit, but not
$m$-limit arrangements} \label{nmv} Consider an irreducible pair
$(\mathcal P ,\mathcal D)$, where $\mathcal D$ is an arrangement of
type $\Gamma_{3,1,1}^{4,0}$ from Theorem \ref{nonreal}. This pair
has ${\rm type}\,(\mathcal P ,\mathcal D) =(m,4,1,4,0)$. In Theorem
\ref{nonreal} it was proven that $\mathcal D$ is not limit for all
$m$. But the line arrangement $\mathcal D$ is virtually $m$-limit
for $m\geq 5$. Indeed, the type of the pair $(\mathcal P ,\mathcal
D)$ has a corresponding data set $(m,4,1,1,0)$. Surfaces $Y$ are of
${\rm type\,}(Y)=(m,4,1,1,0)$, that is, surfaces of degree $m$, the
double curve $D$ of which is an elliptic curve of degree 4, exist.
Since $D$ is a complete intersection of two quadrics $Q_{1}$ and
$Q_{2}$, we can take a surface $Y$, given by equation $F=0$, where
the polynomial $F$ is defined by formula (\ref{Q}). Such surfaces
are completely degenerative, but the limit double curve $\mathcal D$
has a graph not of type $\Gamma_{3,1,1}^{4,0}$, but a graph of type
depicted in Fig. 3.

\subsection{Examples of virtually not limit  line arrangements}
Consider an arrangement of $m$ planes $\mathcal P ={\mathcal
P}_{m}$ in general position in ${\mathbb P}^{3}$. Let ${\mathcal
P}_{m-1} \subset {\mathcal P}_{m}$ be an arrangement of some $m-1$
of these planes, and $\mathcal D ={\mathcal L}_{m-1}$ be the
double curve of the surface ${\mathcal P}_{m-1} $. Then $\bar
d={m-1\choose 2}$, $\tau_0=\tau_2=0$, $\tau_1={m-1\choose 2}$,
$\tau_3={m-1\choose 3}$ and the pair $(\mathcal P,\mathcal D)$ has
type

\begin{equation}\label{Lm}
\text{\rm type}(\mathcal P,\mathcal D)=\left( m,{m-1\choose
2},{m-1\choose 2},0,{m-1\choose 3}\right)
\end{equation}

\begin{prop}\label{nonreal2}
If $(\mathcal P,\mathcal D)$ is an irreducible pair of type {\rm
(\ref{Lm})} and $m\geq 5$, then $\mathcal D$ is not virtually
limit.
\end{prop}
\proof Assume that $\mathcal D$ is virtually limit. Then there
exists a surface $Y$ of type
$${\rm type\,}(Y)=\left( m,{m-1\choose 2},{m-1\choose
2},0,{m-1\choose 3}\right)$$ corresponding to the type of the pair
$(\mathcal P,\mathcal D)$.

Consider a generic projection of the surface $Y$ to the plane.
Then, by Proposition \ref{propx}, the branch curve $B$ of the
generic projection has the following invariants: $\deg B =
2d=2(m-1)$, $g-1 ={m-1\choose 2}-m+1$, $c=3{m-1\choose 2}$, $n=0$.
But then the degree the curve $\check{B}$ dual to the branch curve
$B$ equals
$$\begin{array}{l}\deg \check{B}=\deg B(\deg
B-1)-3c= \\ 2(m-1)(2m-3)-9{m-1\choose 2}=(3-\frac{1}{2}m)(m-1)
\leq 2\end{array}$$ for $m\geq 5$, and it is impossible. \qed

\subsection{Examples of surfaces completely degenerative only virtually}
Consider a degenerative surface $Y$ from section \ref{examp}. The
double curve $D$ of $Y$ is the Hartshorne's curve. The curve $D$
is smooth, has degree $\bar d =30$ and genus $\bar g =113$. Then
$\text{\rm type}(Y)=(m,30,1,113,0)$ and, consequently, $\tau_{2} =
\bar d + \bar g -k =142$ and the collection corresponding to the
type of $Y$ is $(m,30,1,142, 0)$. Let us show that the surface $Y$
is virtually degenerative.

\begin{prop} For $m\geq 31$  there exists an irreducible pairs of
$\text{type}(\mathcal P,\mathcal D)=(m,30,1,142, 0)$.
\end{prop}
\proof Consider a graph $\Gamma(\mathcal D)$ depicted on Fig. 4.

\begin{figure}[h]
\begin{picture}(0,140)(190,-90)
\put(70,13){$\cdot$} \put(70,3){$\cdot$} \put(70,-7){$\cdot$}
\put(70,33){\circle*{3}} \put(57,33){$v_1$} \put(55,-27){$v_{11}$}
\put(70,-27){\circle*{3}} \put(100,3){\circle*{3}}
\put(99,-7){$v_{12}$} \put(70,33){\line(1,-1){30}}
\put(70,-27){\line(1,1){30}} \put(100,3){\line(1,0){35}}
\put(135,3){\circle*{3}} \put(130,-7){$v_{13}$}
\put(135,3){\line(1,0){35}} \put(170,3){\circle*{3}}
\put(165,-7){$v_{14}$} \put(200,-7){$v_{15}$}
\put(170,3){\line(1,0){35}} \put(205,3){\circle*{3}}
\put(205,3){\line(1,0){35}} \put(240,3){\circle*{3}}
\put(240,3){\line(1,0){35}} \put(275,3){\circle*{3}}
\put(305,33){\circle*{3}} \put(305,-27){\circle*{3}}
\put(275,3){\line(1,-1){30}} \put(275,3){\line(1,1){30}}
\put(305,13){$\cdot$} \put(305,3){$\cdot$} \put(305,-7){$\cdot$}
\put(235,-7){$v_{16}$} \put(265,-7){$v_{17}$} \put(310,32){$v_{18}$}
\put(310,-30){$v_{28}$} \put(135,3){\line(0,1){30}}
\put(170,3){\line(0,1){30}} \put(205,3){\line(0,1){30}}
\put(135,33){\circle*{3}} \put(170,33){\circle*{3}}
\put(205,33){\circle*{3}} \put(132,38){$v_{29}$}
\put(167,38){$v_{30}$} \put(202,38){$v_{31}$}
\put(143,-50){$\Gamma^{(30,0)}_{(25,1,3,0,0,0,0,0,0,0,0,2)}$}
\put(173,-80){$\text{\rm Fig}.\, 4$}
\end{picture}
\end{figure}

\noindent It has $31$ vertices, two of which ($v_{12}$ and
$v_{17}$) are of valence $12$, three vertices ($v_{13}$, $v_{14}$,
$v_{15}$) are of valence $3$, the rest vertices are of valence
$1$. Applying (\ref{v}), we find $\tau_{2} =142$. Therefore,
$\Gamma(\mathcal D)$ really is of the mentioned type. It is easy
to see that the pair$(\mathcal P ,\mathcal D)$ is irreducible,
that is, the graph of the complementary curve $\mathcal R=\mathcal
L\setminus \mathcal D$ is connected. \qed

\subsection{Examples of virtually not degenerative surfaces} \label{example2}
Below we prove the existence of a surface $Y\subset {\mathbb P}^3$
with ordinary singularities, whose double curve $D$ has unique
triple point and consists of three components: two lines $L_1$ and
$L_2$ and a conic $Q$. Such a curve $D$ lies on the union of two
planes $P_1$ and $P_2$ such that $L_1\cup L_2\subset P_1$,
$Q\subset P_2$, $s=L_1\cap L_2 \cap Q\in P_1\cap P_2=L$, $L_i\neq
L$ for $i=1,2$ and $Q$ transversally intersect the line $L$.

Let us show that any surface $Y$ with ordinary singularities with
double curve $D$, described above, can not be virtually completely
degenerated. Indeed, the type of such a surface $Y$ is
$\text{type}(Y)=(m,4,3,0,1)$. If $Y$ is virtually completely
degenerated, then an irreducible pair $(\mathcal P,\mathcal D)$
with $\text{type}(\mathcal P,\mathcal D)=(m,4,3,1,1)$ must exist,
that is, there must exist four double lines $\mathcal D$ of a
plane arrangement $\mathcal P$ having the following invariants:
$\tau_2=\tau_3=1$ and, moreover,  after removing the triple point
the curve $\mathcal D$ is decomposed into three ($k=3$) connected
components.

Note that there is a degeneration of the curve $D$ into the union
of four lines having one triple point and one double point (to
obtain such a degeneration we need to degenerate the conic into a
pair of lines lying in $P_2$ so that one of the lines passes
through the point $s$ and the other one does not).

Nevertheless, it is easy to see that there does not exist such
union of double curves of any plane arrangement $\mathcal P=\cup
P_i$ in general position. Indeed, three lines meeting at the
triple point must be pairwise intersections of three planes, say
$P_1$, $P_2$, and $P_3$. It follows from conditions $k=3$ and
$\tau_2=1$ that the fourth line must intersect one (and only one)
of these lines (without loss of generality we can assume that the
fourth line intersects the line $L_{1,2}$). It means that the
fourth line is the intersection of the fourth plane $P_4$ of the
arrangement $\mathcal P$ with one of two planes $P_1$ or $P_2$ and
then this line must intersect either $L_{1,3}$ or $L_{2,3}$. But
it is impossible, since $\tau_2=1$.

Let us prove the existence of a surface $Y$ with ordinary
singularities, $\deg Y\geq 8$, the double curve of which is $D$.
For this let us blow up the point $s$ in $\mathbb P^3$ and after
that let us consecutively blow up the proper transforms of the
curves $L_1$, $L_2$, and $Q$. Let $\sigma :\widetilde{\mathbb
P}^3\to \mathbb P^3$ be a composition of these monoidal
transformations, $E$, $E_1$, $E_2$, $E_3\subset \widetilde{\mathbb
P}^3$ the proper transforms of the exceptional divisors of each of
these blow ups, and let $\widetilde P=\sigma^*(P)$ be the total
inverse image of a plane $P\subset \mathbb P^3$.

Let us show that a generic member $\widetilde X$ of the linear
system $\mid m\widetilde P-3E-2E_1-2E_2-2E_3\mid $ is a smooth
surface if $m\geq 8$. For this let us choose a system of
homogeneous coordinates in $\mathbb P^3$ as follows (remind that
the systems of homogeneous coordinates in $\mathbb P^3$ are
uniquely determined by the choice of ordered sets of four planes
in general position). We take $P_1$ and $P_2$ as the first two
coordinate planes. As the third plane, we take any plane $P_3\neq
P_2$ passing through $s$ and touching the conic $Q$ at this point.
We choose the fourth coordinate plane $P_4$ such that it is in
general position with $P_1$, $P_2$, and $P_3$ and touches the
conic $Q$ at some point $s_1\in Q$, $s_1\neq s$. In the coordinate
system chosen in such a way, the lines $L_i$, $i=1,2$, are given
by equations
$$ z_1=(c_iz_2+z_3)=0, \qquad c_1\neq c_2, $$
and the curve $Q$ is given by
$$ z_2 = (z_3z_4-c_3z_1^2)  = 0, \qquad c_3\neq 0 .$$
The triple point $s$ of $D$ has coordinates $(0:0:0:1)$.

Denote by $F_i$ the following homogeneous polynomials: $F_1=z_1$,
$F_2=z_3z_4+c_1z_2z_4-c_3z_1^2$, $F_3=z_3z_4+c_2z_2z_4-c_3z_1^2$,
$F_4=z_4$ and consider a linear system of surfaces $\overline
S_{a_1,a_2,a_3}\subset\mathbb P^3$ given by homogeneous equation
$$F_1F_2F_3F_4^3+
a_1F_1^2F_2^2F_4^2+a_2F_1^2F_3^2F_4^2+a_3F_2^2F_3^2=0.$$ Denote
also by $x_i=\frac{z_i}{z_4}$ nonhomogeneous coordinates in the
chart $V_4=\mathbb P^3\setminus P_4\simeq \mathbb C^3$ and put
$S_{a_1,a_2,a_3}=\overline S_{a_1,a_2,a_3}\cap V_4$.
\begin{claim} \label{claim} {\rm (i)} For almost all points $(a_1,a_2,a_3)$
the proper transforms \newline
$\sigma^{-1}(S_{a_1,a_2,a_3})\subset \sigma^{-1}(V_4)$ are
nonsingular surfaces.

{\rm (ii)} For $j=1,2,3$ and for almost all points $(a_1,a_2,a_3)$
the intersections \newline $\sigma^{-1}(S_{a_1,a_2,a_3})\cap E_j$
are nonsingular $2$-sections of the ruled surfaces $E_j$, and the
intersections $\sigma^{-1}(s)\cap \sigma^{-1}(S_{a_1,a_2,a_3})\cap
E_j$ consist of two points.

{\rm (iii)} The base locus of the linear system
$\sigma^{-1}(S_{a_1,a_2,a_3})$ consists of three rational curves
lying in $E$, their images after blow down of the divisors $E_1$,
$E_2$, and $E_3$ are lines in the exceptional divisor
$E^{\prime}\simeq \mathbb P^2$ of the blow up of $s$.  In some
analytic neighborhood of $E^{\prime}$, after the blow down the image
of a generic surface of this linear system is decomposed into three
irreducible nonsingular components each of which intersects
transversally $E^{\prime}$ along one of three lines in
$E^{\prime}\simeq \mathbb P^2$ being the image of the base locus of
the linear system.
\end{claim}

\proof In coordinates $x_1$, $x_2$, $x_3$, the linear system of
the surfaces $S_{a_1,a_2,a_3}$ is given by
$$\begin{array}{l} x_1(x_3+c_1x_2-c_3x_1^2)(x_3+c_2x_2-c_3x_1^2)+
a_1x_1^2(x_3+c_1x_2-c_3x_1^2)^2+
\\
a_2x_1^2(x_3+c_2x_2-c_3x_1^2)^2+
a_3(x_3+c_1x_2-c_3x_1^2)^2(x_3+c_2x_2-c_3x_1^2)^2=0.\end{array}$$
Let us take new coordinates in $V_4$:
$$\begin{array}{rl} y_1 & =x_1
\\ y_2 & =x_3+c_1x_2-c_3x_1^2 \\ y_3 & =x_3+c_2x_2-c_3x_1^2.\end{array} $$
In these coordinates the linear system of the surfaces
$S_{a_1,a_2,a_3}$ is given by
$$y_1y_2y_3+a_1y_1^2y_2^2+a_2y_1^2y_3^2+a_3y_2^2y_3^2=0,$$ that is, for almost all
points $(a_1,a_2,a_3)$ the surfaces $S_{a_1,a_2,a_3}\cap V_4$ are
affine parts of the images of Veronese surfaces under a generic
projection into $\mathbb P^3$ for which the claim is well known.
\qed \\

To show that a generic member  $\widetilde X$ of the linear system
$\mid m\widetilde P-3E-2E_1-2E_2-2E_3\mid $ is a nonsingular
surface for $m\geq 8$, let us note that the surfaces $\widetilde
S_{a_1,a_2,a_3}=\sigma^{-1}(\overline S_{a_1,a_2,a_3})$ belong to
the linear system $\mid 8\widetilde P-3E-2E_1-2E_2-2E_3\mid $. It
follows from  Claim \ref{claim} that for any point $p\subset
\mathbb P^3$, $p\neq s$, after the change of the coordinate plane
$P_4$ by another plane not passing through $p$, the base locus of
the linear system $\mid m\widetilde P-3E-2E_1-2E_2-2E_3\mid $ does
not meet the proper transform $\sigma^{-1}(p)$. Consequently, for
any $m\geq 8$ the linear system $\mid m\widetilde
P-3E-2E_1-2E_2-2E_3\mid $ has the same base locus as for $m=8$,
since the linear system $\mid (m-8)\widetilde P\mid $ does not
have the base points. Finally, it follows from Bertini Theorem
that the generic member $\widetilde X$ of $\mid m\widetilde
P-3E-2E_1-2E_2-2E_3\mid $ satisfies the same properties ($i$)
--- ($iii$) of Claim \ref{claim} as the generic member
$\sigma^{-1}(S_{a_1,a_2,a_3})$ has. From this it follows that the
image $Y=\sigma(\widetilde X)$ of the generic member  $\widetilde
X$ is a surface with ordinary singularities the double curve of
which is $D$.

\section{Concluding remarks} In this section, we formulate some
open questions related to the existence problem of complete
degenerations of surfaces with ordinary singular points.

\subsection{Absolute and relative  complete non degeneracy}
Let $X$ be a smooth projective surface and $g:X\to \mathbb P^3$
some "immersion"\, (that is, $Y=g(X)$ is a surface with ordinary
singularities and the morphism $g:X\to Y$ is a normalization of
$Y$). In subsections \ref{example1} and \ref{example2}, we have
given examples of surfaces $Y$ which can not be completely
degenerate. The reason of impossibility to be completely
degenerate can originate from that the "immersion"\, $g$ is
"bad"\, and for some other "immersion"\, of $X$ its image,
nevertheless, can be completely degenerated. It is possible the
second case: for all "immersions"\, of a surface $X$ in $\mathbb
P^3$ it is impossible to completely degenerate its image $Y$, in
other words, the reason of impossibility of complete degenerations
is in topology of $X$. In the first case we say that $X$ is {\it
relatively completely degenerative}, and in the second case $X$ is
{\it absolutely completely non degenerative}. We say also that $X$
is {\it absolutely completely degenerative} if for any
"immersion"\, $g$ its image $g(X)=Y$ is completely degenerative.

\begin{problem} \label{problem2} {\rm ($i$)} Do there exist absolutely completely
non degenerative surfaces $X${\rm ?} \newline {\rm ($ii$)} Do
there exist absolutely completely degenerative surfaces $X${\rm ?}
\end{problem}

In the case of negative answer on any of these problems the
normalizations $X$ of completely non degenerative  and completely
degenerative surfaces $Y$, described in the article, would give
examples  of relatively completely degenerative surfaces.

\subsection{Problem of adjacency}
Denote by $\mathcal H_{\text{type}(Y)}\subset \mathbb
P^{{m+3\choose 3}-1}$ the quasiprojective variety parametrising
the surfaces with ordinary singularities of the same type as the
type of a surface $Y$, $\deg Y=m$. Let $\Pi_m\subset \mathbb
P^{{m+3\choose 3}-1}$ be the variety parametrising arrangements of
$m$ planes in $\mathbb P^3$ in general position, $\dim \Pi_m=3m$.
It follows from complete degeneracy of $Y$ that $\Pi_m$ and the
closure of $\mathcal H_{\text{type}(Y)}$ have nonempty
intersection.
\begin{problem} Let $Y$ be completely degenerative surface with ordinary singularities.
Is it true that $\Pi_m$ lies in the closure of $\mathcal H_{{\rm
\text{type}}(Y)}${\rm ?}
\end{problem}

This is a part of the following more general problem: to describe
the natural stratification (according to the types of double curves)
of the variety  $\mathcal H_m$ of surfaces in $\mathbb P^3$ with
ordinary singularities  and the adjacencies of these strata.

\subsection{Uniqueness complete degeneracy problem}
We say that two pairs $(\mathcal P,\mathcal D_1)$ and $(\mathcal
P,\mathcal D_2)$ are {\it deformation equivalent} if these pairs
are fibres of flat families $\mathcal D_u\subset \mathcal
P_u\subset \mathbb P^3$, $u\in U$, of plane arrangements in
general position and the configurations of double lines contained
in the plane arrangements. It is obvious that pairs $(\mathcal
P_1,\mathcal D_1)$ and $(\mathcal P_2,\mathcal D_2)$ are
deformation equivalent if and only if $\deg \mathcal P_1=\deg
\mathcal P_2$ and the graphs $\Gamma (\mathcal D_1)$ and $\Gamma
(\mathcal D_2)$ are isomorphic.

Let $Y$ be a surface with ordinary singularities. We say that $Y$
has {\it a unique complete degeneration} if it is completely
degenerative and any two its complete degenerations are
deformation equivalent.
\begin{claim} \label{claim2} Any surface $Y$ with ordinary singularities whose double curve
is of degree not more than four possesses not more than unique
complete degeneration
\end{claim}

\begin{figure}[h]
\begin{picture}(0,400)(160,-360)
\put(10,8){\circle*{3}} \put(10,8){\line(1,0){30}}
\put(40,8){\circle*{3}} \put(0,-20){$(m,1,1,0,0)$}
\put(100,18){\circle*{3}} \put(100,8){\circle*{3}}
\put(130,18){\circle*{3}} \put(130,8){\circle*{3}}
\put(100,18){\line(1,0){30}} \put(100,8){\line(1,0){30}}
\put(90,-20){$(m,2,2,0,0)$} \put(180,8){\circle*{3}}
\put(180,8){\line(1,0){30}} \put(210,8){\circle*{3}}
\put(210,8){\line(1,0){30}} \put(240,8){\circle*{3}}
\put(180,-20){$(m,2,1,0,0)$} \put(290,28){\circle*{3}}
\put(290,28){\line(1,0){30}} \put(320,28){\circle*{3}}
\put(290,18){\circle*{3}} \put(290,8){\circle*{3}}
\put(320,18){\circle*{3}} \put(320,8){\circle*{3}}
\put(290,18){\line(1,0){30}} \put(290,8){\line(1,0){30}}
\put(275,-20){$(m,3,3,0,0)$} \put(15,-62){\circle*{3}}
\put(45,-62){\circle*{3}} \put(15,-61){\line(1,0){30}}
\put(-0,-72){\circle*{3}} \put(0,-72){\line(1,0){30}}
\put(30,-72){\circle*{3}} \put(30,-72){\line(1,0){30}}
\put(60,-72){\circle*{3}} \put(0,-100){$(m,3,2,0,0)$}
\put(90,-72){\circle*{3}} \put(90,-72){\line(1,0){25}}
\put(115,-72){\circle*{3}} \put(115,-72){\line(1,0){25}}
\put(140,-72){\circle*{3}} \put(140,-72){\line(1,0){25}}
\put(165,-72){\circle*{3}} \put(100,-100){$(m,3,1,0,0)$}
\put(200,-72){\circle*{3}} \put(200,-72){\line(1,0){40}}
\put(240,-72){\circle*{3}} \put(220,-52){\circle*{3}}
\put(200,-72){\line(1,1){20}} \put(240,-72){\line(-1,1){20}}
\put(190,-100){$(m,3,3,0,1)$} \put(15,-122){\circle*{3}}
\put(15,-121){\line(1,0){30}} \put(45,-122){\circle*{3}}
\put(15,-132){\circle*{3}} \put(15,-132){\line(1,0){30}}
\put(45,-132){\circle*{3}}\put(15,-142){\circle*{3}}
\put(15,-153){\circle*{3}} \put(45,-142){\circle*{3}}
\put(45,-153){\circle*{3}} \put(15,-141){\line(1,0){30}}
\put(15,-152){\line(1,0){30}} \put(0,-180){$(m,4,4,0,0)$}
\put(105,-132){\circle*{3}}  \put(135,-132){\circle*{3}}
\put(105,-132){\line(1,0){30}} \put(105,-142){\circle*{3}}
\put(135,-142){\circle*{3}} \put(105,-141){\line(1,0){30}}
\put(90,-153){\circle*{3}} \put(90,-152){\line(1,0){30}}
\put(120,-153){\circle*{3}} \put(120,-152){\line(1,0){30}}
\put(150,-153){\circle*{3}} \put(90,-180){$(m,4,3,0,0)$}
\put(295,-52){\circle*{3}} \put(295,-72){\circle*{3}}
\put(305,-61){\circle*{3}} \put(295,-52){\line(1,-1){10}}
\put(295,-72){\line(1,1){10}} \put(305,-61){\line(1,0){15}}
\put(320,-61){\circle*{3}} \put(280,-100){$(m,3,1,1,0)$}
\put(180,-142){\circle*{3}} \put(180,-141){\line(1,0){30}}
\put(210,-142){\circle*{3}} \put(210,-141){\line(1,0){30}}
\put(240,-142){\circle*{3}} \put(180,-153){\circle*{3}}
\put(180,-152){\line(1,0){30}} \put(210,-153){\circle*{3}}
\put(210,-152){\line(1,0){30}} \put(240,-153){\circle*{3}}
\put(180,-180){$(m,4,2,0,0)$} \put(300,-142){\circle*{3}}
\put(300,-141){\line(1,0){25}} \put(325,-142){\circle*{3}}
\put(280,-153){\circle*{3}} \put(280,-152){\line(1,0){20}}
\put(300,-153){\circle*{3}} \put(300,-152){\line(1,0){20}}
\put(320,-153){\circle*{3}} \put(320,-152){\line(1,0){20}}
\put(340,-153){\circle*{3}} \put(279,-180){$(m,4,2,0,0)$}
\put(0,-232){\circle*{3}} \put(0,-232){\line(1,0){20}}
\put(20,-232){\circle*{3}} \put(20,-232){\line(1,0){20}}
\put(40,-232){\circle*{3}} \put(40,-232){\line(1,0){20}}
\put(60,-232){\circle*{3}} \put(60,-232){\line(1,0){20}}
\put(80,-232){\circle*{3}} \put(10,-260){$(m,4,1,0,0)$}
\put(160,-217){\circle*{3}} \put(160,-237){\circle*{3}}
\put(170,-227){\circle*{3}} \put(160,-217){\line(1,-1){10}}
\put(160,-237){\line(1,1){10}} \put(170,-227){\line(1,0){15}}
\put(185,-227){\circle*{3}} \put(185,-227){\line(1,0){15}}
\put(200,-227){\circle*{3}} \put(150,-260){$(m,4,1,1,0)$}
\put(215,-297){\line(1,-1){15}} \put(230,-312){\circle*{3}}
\put(200,-282){\circle*{3}} \put(200,-312){\circle*{3}}
\put(215,-297){\circle*{3}} \put(200,-282){\line(1,-1){15}}
\put(200,-312){\line(1,1){15}} \put(215,-296){\line(1,1){15}}
\put(230,-282){\circle*{3}} \put(185,-340){$(m,4,1,3,0)$}
\put(295,-203){\circle*{3}} \put(325,-203){\circle*{3}}
\put(325,-232){\circle*{3}} \put(295,-232){\circle*{3}}
\put(295,-202){\line(1,0){30}} \put(295,-202){\line(0,-1){30}}
\put(295,-232){\line(1,0){30}} \put(325,-202){\line(0,-1){30}}
\put(280,-260){$(m,4,1,1,0)$} \put(80,-282){\line(0,-1){30}}
\put(80,-282){\circle*{3}} \put(80,-312){\circle*{3}}
\put(95,-297){\circle*{3}} \put(80,-282){\line(1,-1){15}}
\put(80,-312){\line(1,1){15}} \put(95,-296){\line(1,0){25}}
\put(120,-297){\circle*{3}} \put(72,-340){$(m,4,2,0,1)$}
\put(152,-370){$\text{\rm Fig}.\, 5$}
\end{picture}
\end{figure}

\proof All possible realizable graphs with the number of edges not
more than four are depicted in Fig. 5. The type $(m,\overline
d,k,\overline g,t)$ of a surface $Y$, in which the plane arrangement
$\mathcal P_m$ in general position can be smoothed outside of the
configuration of corresponding double curve $\mathcal D$, is written
under the each graph.

One can see from this list of graphs that a surface $Y$ could have
more than one complete degeneration only in two cases when
$\text{type}(Y)=(m,4,2,0,0)$ or $\text{type}(Y)=(m,4,1,1,0)$. But in
the first case graphs $\Gamma(\mathcal D)$ have different types
$\Gamma^{(2,0)^2}_{(2,1)^2}$ and $\Gamma^{(1,0)(3,0)}_{(2)(2,2)}$,
and if a plane arrangement $\mathcal P$ is smoothed outside of
$\mathcal D$ with graph $\Gamma^{(2,0)^2}_{(2,1)^2}$, than the
double curve $D$ of $Y$ consists of two irreducible components of
degree two; and if the smoothing takes place outside of  $\mathcal
D$ with graph $\Gamma^{(1,0)(3,0)}_{(2)(2,2)}$, then the irreducible
components of $D$ have degrees one and three, that is, in this case
the regenerated surfaces $Y$ have different (extended) types. In the
second case the graphs $\Gamma(\mathcal D)$ have the same type
$\Gamma^{(4,0)}_{(3,1,1)}$, but plane arrangements  $\mathcal P$ can
not be smoothed outside of one of configurations of double curves
corresponding to this type (see Theorem \ref{nonreal}).

\begin{problem} Does any surface $Y$ possess not more than unique complete degeneration{\rm ?}
\end{problem}

\subsection{Smoothings in symplectic case} Above we gave many examples
of pairs $(\mathcal P,\mathcal D)$ which can not be smoothed
outside of $\mathcal D$. In some cases the obstructions to
smoothing were purely topological (for example, the negativity of
degree of the dual curve of the branch curve $B\subset \mathbb
P^2$ of generic projection to the plane of smoothed surface $Y$),
in the other cases the obstructions, possibly, have algebraic
geometry nature (for example, the pairs $(\mathcal P,\mathcal D)$
with the curve $\mathcal D$ whose graph is depicted in Fig. 4).

To understand better the nature of these obstructions, it is useful
to generalize the problem of complete degenerations of algebraic
surfaces with ordinary singularities to the case of symplectic
varieties. Namely, we say that a compact real four dimensional
subvariety $M$ of $\mathbb C\mathbb P^3$ is a {\it symplectic
variety with ordinary singularities} if for each point $p\in M$
there is a neighborhood $V$ of $p$ such that either variety $V\cap
M$ is decomposed into $n$ smooth components ($n\leq 3$) being
symplectic submanifolds of $\mathbb C\mathbb P^3$ with respect to
the Fubini-Studi symplectic form and meeting transversally along
smooth symplectic surfaces ("double curves"\, of $M$; in the case
$n=3$ the point $p$ is a triple point of $M$) or $V\cap M$ is a
complex analytic variety given in some complex analytic coordinates
in $V$ by equation $x^2-yz^2=0$ (and in this case the point $p$ is
called a pinch of $M$). A definition of complete degeneration of
varieties with ordinary singularities can be also generalized for
symplectic varieties.
\begin{problem} Do there exist pairs $(\mathcal P,\mathcal D)$
such that the plane arrangement $\mathcal P$ can not be smoothed
outside of $\mathcal D$ in the context of algebraic geometry, but
this arrangement can be smoothed outside of $\mathcal D$ in
symplectic context{\rm ?}
\end{problem}

\end{document}